\documentclass[reqno,12pt]{article}
\usepackage{latexsym}
\usepackage{amsmath}
\usepackage{amssymb}
\usepackage{bbm}
\usepackage{amsthm}
\usepackage[authoryear]{natbib}
\usepackage{color, enumerate}
\usepackage{graphicx} 
\usepackage{epsfig}

\usepackage{mathtools}

 0 \DeclareMathAlphabet{\mathbfit}{OT1}{cmr}{bx}{it}

\newcommand{\bea}{\begin{eqnarray*}}
\newcommand{\eea}{\end{eqnarray*}}
\newcommand{\be}{\begin{eqnarray}}
\newcommand{\ee}{\end{eqnarray}}
\newcommand{\beq}{\begin{equation}}
\newcommand{\eeq}{\end{equation}}
\providecommand{\mb}[1]{\mathbb{#1}}

\providecommand{\Abs}[1]{\Bigr\lvert#1\Bigl\rvert}
\providecommand{\norm}[1]{\lVert#1\rVert}

\newcommand{\dd}{\mathrm{d}}

\newtheorem{thm}{Theorem}
\newtheorem{cor}[thm]{Corollary}
\newtheorem{lem}[thm]{Lemma}

\newtheorem{prop}[thm]{Proposition}

\definecolor{Red}{rgb}{1,0,0}

\definecolor{Blue}{rgb}{0,0,1}

\begin{document}

\title{Einstein relation and steady states for the random conductance model} 

\author{
{\small Nina Gantert}\\
\and
{\small Xiaoqin Guo}\\
\and
{\small Jan Nagel}\\
}

\maketitle


Abstract: We consider random walk among iid, uniformly elliptic conductances on $\mathbb{Z}^d$, and prove the Einstein relation (see Theorem \ref{ER}). It says that the derivative of the velocity of a biased walk as a function of the bias equals the diffusivity in equilibrium.
For fixed bias, we show that there is an invariant measure for the environment seen from the particle. These invariant measures are often called steady states.
The Einstein relation follows at least for $d \geq 3$, from an expansion of the steady states as a function of the bias (see Theorem~\ref{newER}), which can be considered our main result. This expansion is proved for $d \geq 3$.
In contrast to \cite{guo2012}, we need not only convergence of the steady states, but an estimate on the rate of convergence (see Theorem \ref{newmain2}).


\section{Introduction}
We consider random walk among iid, uniformly elliptic random conductances.
More precisely, let
$\mathcal{B}(\mathbb{Z}^d)$ be the set of non-oriented nearest-neighbor bonds in $\mathbb{Z}^d$, $d \geq 2$ and $\Omega:=(0,\infty)^{\mathcal{B}(\mathbb{Z}^d)}$. 
An element $\omega\in\Omega$ is called an {\it environment}. For $(x,y) \in \mathcal{B}(\mathbb{Z}^d)$, the weight $\omega(x,y)$ is the {\sl conductance} of the bond $(x,y)$.
For $\omega \in \Omega$, the random walk in the environment $\omega$ is the Markov process $(X_n)_{n\geq 1}$ with transition probabilities 
\begin{align*}
P^x_\omega(X_{n+1}=y+e | X_n=y) = \frac{\omega(y,y+e)}{\sum_{|e'|=1} \omega(y,y+e')}
\end{align*}
and $P^x_\omega (X_0=x)=1$. In other words, the transition probabilities are proportional to the conductances of the bonds.
The distribution $P_\omega^x$ of the random walk is called \emph{quenched law}. For a probability measure $P$ on $\Omega$, the averaged measure $\mathbb{P}^x= P \times P_\omega^x: = \int P_\omega^x \, P(d\omega) $ is called the \emph{annealed law}. We will assume throughout this paper that
\begin{itemize}
\item[(i)] there is a $\kappa>1$ such that $ \kappa^{-1}\leq \omega(x,y) \leq \kappa \text{ for all } (x,y) \in \mathcal{B}(\mathbb{Z}^d)$,
\item[(ii)] $(\omega(e))_{e\in\mathcal B(\mathbb Z^d)}$ are iid under $P$.
\end{itemize}
Now, for $\lambda \in (0,1)$ and $\ell \in \mathbb{R}^d$, $|\ell|=1$, define the perturbed environment $\omega^\lambda$ of $\omega \in \Omega$ by
\begin{align*}
\omega^\lambda (x,y) = \omega(x,y)e^{\lambda \ell \cdot (x+y)} ,
\end{align*} 
where $``\cdot"$ denotes the scalar product in $\mathbb R ^d$.
We denote the quenched and annealed measures of the random walks in the perturbed environment as $P^x_{\omega,\lambda}:= P^x_{\omega^\lambda}$ and $\mathbb{P}_\lambda^x:= P \times P_{\omega,\lambda}^x$, respectively.
If the starting point $x$ is the origin, we will omit the superscript $x=0$, for instance, we write $P_\omega$ instead of $P^0_\omega$ and $\mathbb{P}$ instead of $\mathbb{P}^0$.
It goes back to \cite{masi1988} that
for $\lambda =0$, the random walk satisfies a functional central limit theorem under $\mathbb{P}$: it converges under diffusive scaling to a Brownian motion $(B_t)_{t\geq 0}$ with a deterministic, diagonal covariance matrix $\Sigma$.
For $\lambda > 0$, it was shown in \cite{shen2002} that the random walk is ballistic, i.e. there is a deterministic vector $v(\lambda)$ with $v(\lambda) \cdot \ell > 0$ such that
\begin{align*}
\lim\limits_{n \to \infty}\frac{X_n}{n} = v(\lambda) \quad  \mathbb{P}_\lambda-\rm{a.s.}
\end{align*}
The Einstein relation says that the derivative of the speed with respect to the perturbation relates to the covariance matrix of the unperturbed random walk among random conductances as follows.

\begin{thm}[Einstein relation] \label{ER}
\begin{align}\label{ERconv}
\lim_{\lambda \to 0} \frac{v(\lambda)}{\lambda} = \Sigma \ell\, .
\end{align}
\end{thm}

We remark that Theorem \ref{ER} holds true for $d=1$ as well, but in that case, it can be shown by explicit calculation.
Theorem \ref{ER} was known for the case when the conductances take only two values and dimension is at least $3$, see \cite{komorowskiolla2005einstein}.
The Einstein relation is conjectured to be true in general for reversible motions which behave diffusively. We refer to 
\cite{einstein1956} for a historical reference and to \cite{spohn} for further explanations.
A weaker form of the Einstein relation holds indeed true under such general assumptions and goes back to \cite{lebrost1994}.
However, \eqref{ERconv} was only established in examples: 
\cite{Loulakis} and \cite{loulakis05} consider a tagged particle in an exclusion process, \cite{komorowskiolla2005mobility} and \cite{lucarenatoflorian} 
investigate other examples of space-time environments,
the paper \cite{komorowskiolla2005einstein} mentioned above gives the result for particular random walks among random conductances, \cite{ganmatpia2012} treats reversible diffusions in random environments and \cite{ahoz2013einstein} considers biased random walks on Galton-Watson trees. The only result, to our best knowledge, for a non-reversible situation is given in \cite{guo2012}. 
For results on the steady states in the case of diffusions, we also refer to  
forthcoming work of Pierre Mathieu and 
Andrey Piatnitski.

Note that in \cite{guo2012}, when the random walk is a martingale, 
the Einstein relation \eqref{ERconv} is a consequence of a more general convergence theorem for the steady states. In the random conductance model, due to the presence of the corrector, to generalize \eqref{ERconv} we need finer estimates of the rate of the convergence of the steady states, for which we first introduce some notation.

When $(\theta_z\omega)(x,y)=\omega(z+x,z+y)$ denotes the environment shifted by $z$, we set $\bar\omega_n:=\theta_{X_n}\omega, n\ge 0$. The Markov chain $(\bar\omega_n)$ is called the \emph{environment seen from the particle} and has generator $L$ acting on bounded functions $f:\Omega \rightarrow \mathbb{R}$ as 
\begin{align*}
(Lf)(\omega) = \sum_{y:|y|=1} P_\omega(X_1=y) (f(\theta_y\omega)-f(\omega)) .
\end{align*}
For $\lambda=0$, the Markov chain $(\bar\omega_n)$ has an invariant measure $Q_0$ given by 
\begin{equation}\label{e16}
\frac{\dd Q_0}{\dd P}(\omega) = Z^{-1} \sum_{y:|y|=1} \omega(0,y) \in [1/\kappa^2, \kappa^2],
\end{equation}
with a normalization constant $Z$. If $\lambda>0$ and the law of $X_n$ is given by $P_{\omega,\lambda}$, an invariant measure $Q_\lambda$ for the Markov chain $\bar\omega_n$ can be defined as the $\mathbb P_\lambda$-a.s. limit 
\begin{align} \label{invmeasurelambda}
Q_\lambda f = \lim_{n\to \infty} \frac{1}{n} \sum_{k=1}^n f(\bar{\omega}_k),
\end{align} 
which has an expression in terms of the regeneration times defined in Section \ref{sec:reg},
see \eqref{e5}.
The invariant measures $Q_\lambda$ are often called ``steady states''.
The Einstein relation \eqref{ERconv} is a consequence of a first order expansion of ${Q}_\lambda$ around $\lambda=0$, see Theorem \ref{newER} below. 

To describe the limit, we let $\mathcal H_{-1}$ denote the set of all  functions $f:\Omega\to\mathbb R$ in  $L^2(Q_0)$ such that the limit
\begin{align*}
\lim_{n\to\infty}\frac{1}{n} Q_0 E_\omega\left[ \left(  \sum_{k=0}^{n} f(\bar\omega_k) \right)^2 \right]=:\sigma^2(f)
\end{align*}
exists and is finite. For a variational characterization of the space $\mathcal{H}_{-1}$, we refer to \cite{kipnis1986central,klo2012fluctuations}. 
In the classical paper \cite{kipnis1986central}, it is proved that for $f\in \mathcal{H}_{-1}$,  the process
\begin{align*}
\left(
\frac{1}{\sqrt{n}} \sum_{k=1}^{nt} f(\bar{\omega}_k)
\right)_{t\ge 0}
\end{align*}
converges weakly (under $\mathbb P_0$) to a Brownian motion $N^f_t$ with variance $\sigma^2(f)$. 
For our result, we  consider the subspace $\mathcal{F}$ of bounded continuous functions $f:\Omega \to \mathbb{R}$ depending only on a finite set of conductances, that is, $f(\omega) = \tilde f((\omega_e)_{e\in E})$ for a finite set $E\subset \mathcal{B}(\mathbb{Z}^d)$. We remark that it follows from \cite{mourrat2011variance} that if $f\in \mathcal{F}$ and $d\geq3$, then $f-Q_0f \in \mathcal{H}_{-1}$. 
Consider the 2-dimensional process 
\begin{align*}
\frac{1}{\sqrt{n}}\left( \sum_{k=1}^n f(\bar{\omega}_k)-Q_0f, \sum_{k=1}^n d(\bar{\omega}_k,0)\cdot \ell \right),
\end{align*}
where $d (\omega,x)=  E_\omega^{x}[(X_1-X_0)]$ is the local drift. If $f-Q_0f\in \mathcal{H}_{-1}$, this process converges by \cite{kipnis1986central} in distribution under $\mathbb{P}_0$ to a 2-dimensional normal random variable ${(N_1^f,N_1^d\cdot\ell)}$. Define then 
\begin{equation}\label{defLambda}
\Lambda f = -{\rm Cov}{ (N_1^f,N_1^d\cdot\ell)}.
\end{equation}

 The following theorem is our main result.
\begin{thm} \label{newER}
If $d\ge 3$, we have
\begin{equation}\label{1stOrderExp}
\lim_{\lambda \to 0} \frac{{Q}_\lambda f - {Q}_0 f}{\lambda}
= \Lambda f
\end{equation}
for any $f\in \mathcal{F}$.
\end{thm}

We remark that it follows by similar arguments as in \cite{guo2012} that $Q_\lambda\Rightarrow Q_0$ for $d\ge 2$. The first order expansion \eqref{1stOrderExp} of the measure $Q_\lambda$ is obtained in \cite{guo2012} for $P$ that satisfies some ballisticity condition, where a regeneration structure creates enough decorrelation for the environments along the path. In our case where the unperturbed environment $P$ is not ballistic, the first order expansion \eqref{1stOrderExp} is more delicate. When $d\ge 3$, making use of the optimal variance decay for 
the environment seen from the particle in \cite[Theorem 1.1]{debuyermourrat2014} (for the unperturbed environment) and a 1-dependent regeneration structure (for perturbed environments), we obtain \eqref{1stOrderExp}. For $d\le 2$,  it is not clear to us whether \eqref{1stOrderExp} still holds, since the environment seen from the particle process decorrelates at a slower rate and our argument does not go through.

In order to prove Theorem \ref{newER}, we will show the following two theorems.
To simplify notation, we will write $X_t$ for $X_{\lfloor t\rfloor}$ and similarly for summation limits or other indices which are 
defined for integer values. In Theorem \ref{newmain1} we center by $Q_0f$ to point out the relation to Theorem \ref{newER}, but we remark that $f\in \mathcal{H}_{-1}$ actually implies $Q_0f=0$.

\medskip

\begin{thm}\label{newmain1}
For any $t\geq 1$ and $f\in \mathcal{H}_{-1}$, we have 
\begin{equation}\label{newmain1f}
\lim_{\lambda \to 0} \frac{ \frac{\lambda^2}{t} \mathbb{E}_{Q,\lambda}  \sum_{k=0}^{t/\lambda^2} f(\bar\omega_k) - Q_0 f}{\lambda}  = \Lambda f ,
\end{equation}
where $\mathbb{E}_{Q,\lambda}$ is the expectation with respect to $Q_0\times P_{\omega,\lambda}$.
\end{thm}

\medskip

\begin{thm}\label{newmain2}
If $d\ge 3$, then for any $f\in \mathcal{F}$ there exist constants $C=C(\kappa,d,f)$ and $\lambda_0 = \lambda_0(\kappa,d)> 0$ such that  for any $t>1$ and $\lambda< \lambda_0$, 
\begin{equation}\label{newmain2f}
\Abs{
\frac{ \frac{\lambda^2}{t} \mathbb{E}_{Q,\lambda} \sum_{k=0}^{t/\lambda^2} f(\bar\omega_k) - Q_\lambda f}{\lambda}
}  
\leq 
\frac{C}{t^{1/4}}.
\end{equation}
\end{thm}

\medskip

The paper is organized as follows: In Section \ref{A-priori estimates}, we collect some a priori estimates whose proofs are deferred to Section \ref{Pape}. We then define, in Section \ref{sec:reg}, a regeneration structure which will enable us to prove Theorem \ref{newmain2}. As in \cite{ganmatpia2012}, we here have to take into account how the regeneration times and regeneration distances depend on $\lambda$.
In Section \ref{sec:gir}, we prove Theorem~\ref{newmain1}: the main ingredient is Girsanov transform and the technique of proof is similar to \cite{ganmatpia2012}. 
In Section \ref{auLLN}, using the regeneration structure, we prove \eqref{newmain2f} for a class of functions that satisfies a nice inequality, namely \eqref{e20}. In Section \ref{Proof of the Einstein relation}, using the variance estimate of \cite{debuyermourrat2014, mourrat2011variance}, 
we prove Theorem~\ref{newmain2} by verifying \eqref{e20} for all $f\in \mathcal F$ when $d\ge 3$.
 We also show how the Einstein relation \eqref{ERconv} follows, at least in the case $d \geq 3$, from Theorem \ref{newER} and give a different argument to show \eqref{ERconv} for any $d$. 
Finally, in Section \ref{Pape}, we prove the estimates listed in Section \ref{A-priori estimates}.



\section{A-priori estimates}\label{A-priori estimates}

Without loss of generality, assume throughout that $\ell\cdot e_1 = \max_i\{\ell\cdot e_i,-\ell\cdot e_i\}$, where $\{e_i: i=1,\ldots, d\}$ denotes the natural basis for $\mathbb Z^d$. Let $\lambda_1 = \lfloor 1/\lambda \rfloor^{-1}$, such that $1/\lambda_1 \in \mathbb{N}$ and $0\leq \lambda_1-\lambda = o(\lambda)$ as $\lambda\to 0$. For $m \in \mathbb{Z}$ and $L\in \mathbb{N}$ we define the 
hyperplane
\begin{align*}
H_{m,L} =\{x \in \mathbb{Z} ^d |\, x \cdot e_1  =mL/\lambda_1 \} 
\end{align*}
and the hitting time
\begin{align*}
T_{m,L} =\inf\{n\geq 0|\, ( X_n-X_0)\in H_{m,L}\} .
\end{align*}
For $x\in \mathbb{R}^d$, we denote by $|x|$ the 1-norm of $x$.

The constants appearing in this paper will be allowed to depend on the dimension $d$ and the ellipticity constant $\kappa$ but we emphasize that they do not depend on $\lambda$. In the proofs, we use $c,C$ to denote generic positive constants whose values may change from line to line.

\begin{lem}[Bounding the probability to go left before going right]\label{hitting}
There exist $L_0 \in \mathbb{N}, \lambda_0 >0$, depending only 
on the ellipticity constant $\kappa$ and the dimension $d$, such that 
\begin{align*}
P_{\omega,\lambda} (T_{1,L} < T_{-1,L}) \geq \frac{2}{3}
\end{align*}
for all $L\geq L_0, 0<\lambda \leq \lambda_0$ and for all $\omega$.
\end{lem}

\medskip

In the following, fix $L_0$ as in Lemma \ref{hitting} and we write shorter $H_m$ for $H_{m,4L_0}$ and $T_m$ for $T_{m,4L_0}$. For the distance between these hyperplanes we write $L_1=4L_0/\lambda_1$, tacitly ignoring the dependence on $\lambda$. 

\medskip

\begin{cor}
[Bounding the probability to go far to the left] 
\label{cor1} 
Let $\lambda_0$ be the same constant as in Lemma~\ref{hitting}.
For any $\omega \in \Omega, n\in\mb N$ and $\lambda\in(0,\lambda_0)$, we have
\[
P_{\omega,\lambda}(T_{-n/4}<\infty)\le 2^{-n}.
\]
\end{cor}

\begin{lem}
[Hitting times of hyperplanes to the right are small with high probability]
\label{hitting2} 
There exist positive constants $c_1,C_1$ such that for any  $\omega$, $n\geq 1$ and $\lambda \in (0,\lambda_0)$ we have
\begin{align*}
P_{\omega,\lambda}\left(T_{n} \geq \tfrac{C_1n}{\lambda^2} \right)\leq  e^{-c_1n} .
\end{align*}
\end{lem}

\begin{lem}
[Bounds on the moments of the maximum of the walk]
\label{maxbound}
For any $p\geq 1$ there exists a positive constant $C_2$, such that for all $\lambda\in (0,\lambda_0)$ and $t\geq 1$ and for any $\omega$,
\begin{align*}
E_{\omega,\lambda}\left[ \max_{0\leq s\leq t} |\lambda X_{s/\lambda^2}|^p \right] \leq C_2 t^p .
\end{align*}
\end{lem}

Corollary~\ref{cor1} is an immediate consequence of Lemma~\ref{hitting}. The proofs of Lemmas \ref{hitting}, \ref{hitting2} and \ref{maxbound} are given in Section~\ref{Pape}.
The following parabolic Harnack inequality will be used several times in our paper. Let $B_R=\{x\in\mathbb Z^d: |x|_2\le R\}$ and $B_R(x):=x+B_R$.
\begin{thm}[Parabolic Harnack inequality]
Fix $R\ge \sqrt d$. Let $a\in\Omega$ be a configuration of conductances such that $\frac{a(e)}{a(e')}\le C_V$ for any two bonds $e,e'$ in $B_{2R}$. Assume that $u:\mathbb Z^d\times\mathbb Z\to\mathbb R_+$ is a nonnegative function that satisfies the parabolic equation
\[\tag{PE}
u(x,n+1)=\sum_y P_a^x(X_1=y)u(y,n) \qquad\text{ in }B_{2R}\times[0, 4R^2+1].
\]
Then there exists a constant $C=C(d,C_V)$ such that
\[\tag{PHI}
\max_{B_R\times[R,2R^2]}u\le C\min_{(x,n)\in B_R\times[3R^2, 4R^2]}[u(x,n)+u(x,n+1)].
\]
\end{thm}
We remark that the bound $a(e)/a(e')<C_V$ implies the {\it volume-doubling condition} (c.f. \cite[Def. 1.1]{delmotte1999}) and uniform ellipticity. 
For our exponentially growing conductances $\omega^\lambda$, we can choose the constant $C_V$ independent of $\lambda\in[0,1]$, as long as we consider $u$ defined on a subgraph of size $C/\lambda$. 
For the conductances with $a(x,x)>c>0$, the theorem is \cite[Theorem 1.7]{delmotte1999}. In our case $a(x,x)=0$, (PHI) is obtained by applying \cite[Theorem 1.7]{delmotte1999} to even-step jumps of the random walk, see the remarks below \cite[Def. 1.3]{delmotte1999}. The above version of (PHI) can be found in \cite[Def. 2.2]{grigor2002harnack}. Note that when $u$ is not a function of time, i.e, $u(x,n)=u(x,m)=:u(x)$ for all $n,m\in\mathbb Z$, then $u$ satisfies the {\it elliptic equation}
\[\tag{EE}
u(x)=\sum_y P_a^x(X_1=y)u(y)
\]
and (PHI) becomes the {\it elliptic Harnack inequality} 
\[\tag{EHI}
\max_{B_R}u\le C\min_{B_R} u.
\]


\section{Regeneration structure}\label{sec:reg}
In this section we will construct a $1$-dependent regeneration structure on the random walk path so that the inter-regeneration distances and inter-regeneration times are roughly of order $1/\lambda$ and $1/\lambda^2$, respectively. The regeneration structure will then imply the estimate of Theorem \ref{newmain2}. For this we fix the function $f\in \mathcal{F}$ and allow the constants in this chapter to depend on $f$.


\subsection{Auxiliary estimates}
\begin{lem}[Transversal fluctuations are not too large]\label{lem-reg1}
There exists a constant $C_3$ so that for all $\lambda\in(0,\lambda_0)$
\[
P_{\omega,\lambda}(|X_{T_{3/4}}|\ge C_3/\lambda)\le 1-C_3^{-1} .
\]
\end{lem}

\proof
By Lemma~\ref{hitting2} and Lemma~\ref{maxbound},  for any $\theta\geq 1$,
\begin{align*}
&\quad P_{\omega,\lambda}(|X_{T_{3/4}}|\geq \theta/\lambda)\\
&\leq 
P_{\omega,\lambda}(T_{3/4}\geq C_1/\lambda^2)+P_{\omega,\lambda}(\max_{0\le s\le C_1/\lambda^2}|X_s|\ge \theta/\lambda)\\
&\le 
e^{-c_1}+C_1C_2/\theta.
\end{align*}
Taking $\theta$ sufficiently large, the lemma follows.\qed

\begin{lem}
[Bounding the exit measure on a hyperplane to the right] 
\label{hypermeasure} 
There exists a probability measure $\mu_{\omega,\lambda,1}$ on $H_{1}$, which is independent of $\sigma\big(\omega(x,y):x\cdot e_1\leq L_0/\lambda_1, y\in\mathbb Z^d\big)$, and a constant $c_4> 0$ such that 
\[
P_{\omega,\lambda}(X_{T_{1}}=\cdot)\ge c_4 \mu_{\omega,\lambda,1}(\cdot)\, .
\]

\end{lem}
\proof
We will prove the lemma by showing that for any $w\in H_{1}$ and $x= e_1 \cdot  3L_0/4\lambda_1$,
\[
P_{\omega,\lambda}(X_{T_{1}}=w)
\geq
C P_{\omega,\lambda}^x (X_{T_{1/4}}=w|T_{1/4}<T_{-1/4}).
\] 

Indeed, for any $w\in H_{1}$,
\begin{align*}
P_{\omega,\lambda}(X_{T_{1}}=w)
&\geq 
\sum_{y:|y-x|< C_3/\lambda} P_{\omega,\lambda}(X_{T_{3/4}}=y)P_{\omega,\lambda}^y(X_{T_{1/4}}=w)\\
&\geq C
\sum_{y:|y-x|< C_3/\lambda} P_{\omega,\lambda}(X_{T_{3/4}} =y)P_{\omega,\lambda}^{x}(X_{T_{1/4}}=w)\\
&\stackrel{\mathclap{Lem.~\ref{lem-reg1}}}{\geq}\  CC_3^{-1} P_{\omega,\lambda}^{x}(X_{T_{1/4}}=w)\\
&\stackrel{\mathclap{Lem.~\ref{hitting}}}{\geq}\ 
CP_{\omega,\lambda}^{x}(X_{T_{1/4}}=w|T_{1/4}<T_{-1/4}),
\end{align*}
where in the second inequality we 
applied the elliptic Harnack inequality (EHI) to the function $u(y):=P^y_{\omega,\lambda}(X_{T_{1/4}}=w)$ in the ball $B_{2C_3/\lambda}(x)$.

\qed


\subsection{Construction of the regeneration time}

Classical regeneration times are usually defined as times when the random walker crosses a certain hyperplane for the first time and then never comes back. To keep the regeneration times robust for small bias, we allow the path to backtrack a distance of order $1/\lambda$. To decouple the trajectory between these regeneration times, we will then use the ``coin trick'' by \cite{cometszeitouni2004}.

The starting point is that by Lemma~\ref{hypermeasure}, the hitting probability $P_{\omega,\lambda}^{x}(X_{T_1}=\cdot)$ of the next hyperplane dominates $c_4$ times a probability measure $\mu^x_{\omega,\lambda,1}$, 
which is independent of the environment  to the left of the hyperplane 
$H^x_0:=x+H_0$
. Hence for $\beta\in(0,c_4)$ the hitting probability can be decomposed as
\begin{align*}
P_{\omega,\lambda}^x(X_{T_{1}}=\cdot)
=\beta\mu_{\omega,\lambda,1}^x(\cdot)+(1-\beta)\mu_{\omega,\lambda,0}^x(\cdot),
\end{align*}
where 
\begin{align*}
\mu_{\omega,\lambda,0}^x(\cdot)
:=
\frac{P_{\omega,\lambda}^x(X_{T_{1}}=\cdot)-\beta\mu_{\omega,\lambda,1}^x(\cdot)}{1-\beta}.
\end{align*}
By Lemma \ref{hypermeasure}, both $\mu_{\omega,\lambda,1}^x$ and $\mu_{\omega,\lambda,0}^x$ are probability measures on 
$H^x_{1}   = \{y \in \mathbb{Z} ^d | (y-x) \cdot e_1  =L_1 \}$.
Let $(\varepsilon_i)_{i=0}^\infty\in\{0,1\}^{\mathbb{N}_0}$ be iid Bernoulli random variables with law $q_\beta$:
\begin{align*}
q_\beta(\varepsilon_i=1)=\beta \text{ and } q_\beta(\varepsilon_i=0)=1-\beta.
\end{align*}
Intuitively, when $X_n$ is at $x\in H_{i}$ the coin $\varepsilon_i$ will determine whether the hitting point of the next hyperplane $H_{i+1}$ is sampled via $\mu_{\omega,\lambda,0}^x$ or $\mu_{\omega,\lambda,1}^x$. Until reaching $H_{i+1}$, the law of the path will then be the original quenched law, conditioned on the predetermined hitting point.

We now give the formal definition of the regeneration times, for which we first define inductively a path measure given a set of ``hitting rules'' as described above. 
Sample the sequence $\varepsilon:=(\varepsilon_i)_{i=0}^\infty$ according to the product measure $q_\beta$ and fix it. Then, define $P_{\omega,\lambda,\varepsilon}^x$
on the paths by the following steps:
\begin{itemize}
\item {\it Step 1.} For $x\in\mathbb{Z}^d$, set
\begin{align*}
P_{\omega,\lambda,\varepsilon}^x(X_0=x)=1.
\end{align*}
and for any $\mathcal{O}\in\sigma(X_1,X_2,\ldots, X_{T_{1}})$, put
\begin{align}\label{e6}
\nu_{\omega,\lambda,\varepsilon_i}^x(\mathcal{O})
:=
\sum_y 
\big[\varepsilon_i\mu_{\omega,\lambda,1}^x(y)+(1-\varepsilon_i)\mu_{\omega,\lambda,0}^x(y)\big]
P_{\omega,\lambda}^x(\mathcal{O}|X_{T_1}=y).
\end{align}
\item {\it Step 2.}
Suppose the $P_{\omega,\lambda,\varepsilon}^x$-law for paths of length$\,\leq n$ is defined. For any path $(x_i)_{i=0}^{n+1}$ with $x_0=x$, define
\begin{align*}
&P_{\omega,\lambda,\varepsilon}^{x}
(X_{n+1}=x_{n+1},\ldots, X_{0}=x_0)\\
&:=
P_{\omega,\lambda,\varepsilon}^{x}(X_I=x_I,\ldots, X_0=x_0)
\nu_{\omega,\lambda,\varepsilon_J}^{x_I}(X_{n+1-I}=x_{n+1},\ldots, X_1=x_{I+1}),
\end{align*}
where
\[
J=\max\{j\ge 0: H_{j}^{x_0}\cap\{x_i, 0\le i\le n\}\neq\emptyset\}
\] 
is the rightmost hyperplane visited by $(x_i)_{i=0}^{n}$ and
\[
I=\min\{0\le i\le n: x_i\in H_{J}^{x_0}\}
\]
is the hitting time to the $J$-th level.
\item {\it Step 3.}
By induction, the law $P_{\omega,\lambda,\varepsilon}^x$
is well-defined for paths of all lengths.
\end{itemize}

Intuitively, whenever the walker visits new hyperplanes $H_i, i\ge 0$, 
we make him flip a coin $\varepsilon_i$. 
If $\varepsilon_i=0$ (or $1$), he then walks
following the law $\nu_{\omega,\lambda,0}$ (or $\nu_{\omega,\lambda,1}$) until he reaches the $(i+1)$-th hyperplane. The regeneration time $\tau_1$ is defined to be the first time of visiting a new hyperplane
$H_k$ such that the outcome $\varepsilon_{k-1}$ of the previous 
coin-tossing is ``$1$" and the path will never backtrack to level 
$H_{k-1}$ in the 
future.

Note that a path sampled by $P_{\omega,\lambda,\varepsilon}^x$ is not a Markov chain, but
the law of $X_\cdot$ under 
\[
\bar{P}_{\omega,\lambda}^x
:=q_\beta\times P_{\omega,\lambda,\varepsilon}^x
\]
coincides with $P_{\omega,\lambda}^x$. That is, 
\begin{equation}\label{e47}
\bar{P}_{\omega,\lambda}^x(X_\cdot\in\cdot)
=
P_{\omega, \lambda}^x(X_\cdot\in\cdot).
\end{equation}
We denote by 
\begin{equation}\label{notPlambdabar}
\bar{\mathbb P}_\lambda
:= P \times\bar{P}_{\omega,\lambda}
\end{equation}
the law of $X_\cdot$ averaged over the coins and the environment. Expectations with respect to $\bar{P}_{\omega,\lambda}^x$ and $\bar{\mathbb P}_\lambda$ are denoted by 
$\bar{E}_{\omega,\lambda}^x$ and $\bar{\mathbb E}_\lambda$, respectively.
Next, for a path $(X_n)_{n\ge 0}$ sampled according to $P_{\omega,\lambda,\varepsilon}$, we will define the regeneration times. 
To be specific, put $S_0=0, M_0=0$,
and define inductively the times $S_k$ and $R_k$ and the distances $M_k$ by
\begin{align} \label{defSRM}
&S_{k+1}=\inf\{T_{n+1}: nL_1\geq M_k \text{ and }\varepsilon_n=1\},\notag \\
&R_{k+1}=S_{k+1}+T_{-1/4}\circ\theta_{S_{k+1}},\\
&M_{k+1}=X_{S_{k+1}}\cdot e_1+N\circ \theta_{S_{k+1}}\cdot L_1, \qquad k\ge 0.\notag
\end{align}
Here $\theta_n$ denotes the time shift of the path, i.e, $\theta_n X_\cdot=(X_{n+i})_{i=0}^\infty$, 
and
\begin{equation}\label{e46}
N:=\inf\{n: nL_1>(X_i-X_0)\cdot e_1 \text{ for all }i\le T_{-1/4}\}.
\end{equation}
Set
\begin{equation}\label{Kdef}
K:=\inf\{k\ge 1: S_k<\infty, R_k=\infty\}\, ,
\end{equation}
\begin{equation}
\tau_1:=S_K \quad \text{ and }\tau_{k+1}=\tau_k+\tau_1\circ\theta_{\tau_k}.
\end{equation}
The times $(\tau_k)_{k\ge 1}$ are called \textit{regeneration times}.

\begin{figure}[h]
\begin{centering}
\includegraphics[width=0.6\textwidth]{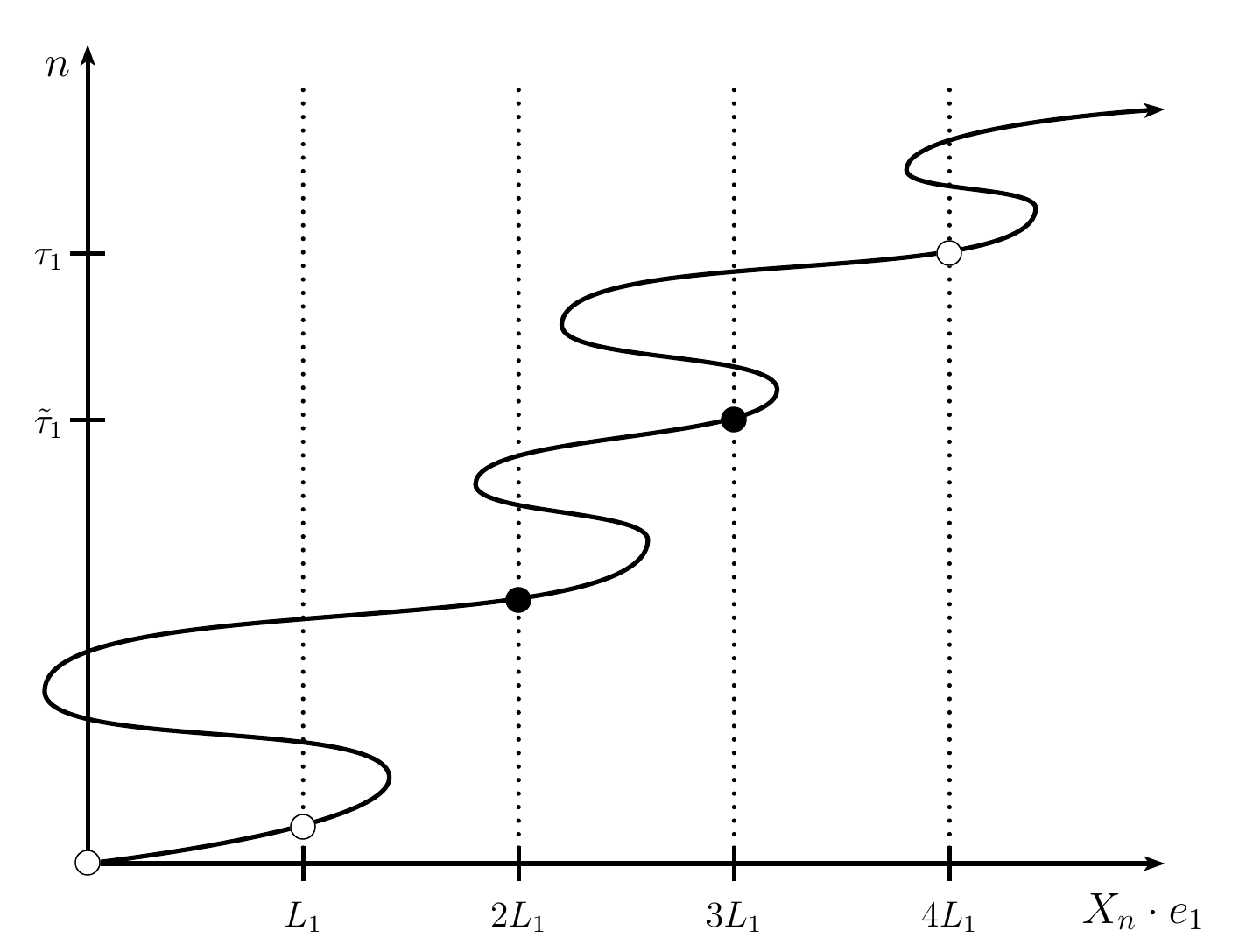}
\caption{A sample path of the random walk. A black dot at $X_{T_n}$ represents a successful coin toss $\varepsilon_n=1$, a white dot corresponds to $\varepsilon_n=0$. After hitting $H_2$, the random walk does not backtrack more than $L_1/4$, but since $\varepsilon_1=0$, this is not a regeneration time. Since $\varepsilon_3=1$ and after reaching $H_4$, the path does not backtrack more than $L_1/4$, we have $\tilde{\tau_1}=T_3$, $\tau_1=T_4$. Note that $\varepsilon_3=1$ implies that the hitting point $X_{T_4}$ was chosen according to $\mu^{X_{T_3}}_{\omega,\lambda,1}$. }
\end{centering}
\end{figure}


\subsection{Renewal property of the regenerations}

The regeneration times possess good renewal properties:
\begin{enumerate}
\item 
Set $\tau_0=0$. For $k\ge 0$, define
\begin{align*}
\tilde{S}_{k+1}:=\inf\{T_{n}: nL_1 \geq M_k \text{ and }\varepsilon_n=1\},
\end{align*}
with $M_k$ is as in \eqref{defSRM}. That is, if we divide the space by the hyperplanes $H_k$ at distance $L_1$ of each other, $\tilde{S}_{k+1}$ is the hitting time of 
a hyperplane
 after 
 the previous maximum $M_k$ 
 is achieved
 and when the coin corresponding to this hyperplane lands head. Note that ${S}_{k+1}$ is the hitting time of the next hyperplane after $\tilde{S}_{k+1}$. Again, with $K$ as in \eqref{Kdef}, set
\begin{align*}
\tilde{\tau}_1:=\tilde{S}_K,\qquad \tilde{\tau}_{k+1}:=\tau_k+\tilde{\tau}_1\circ\theta_{\tau_k}\quad  (k>1).
\end{align*}
Namely, $\tilde\tau_k$ is the hitting time to the previous hyperplane of $X_{\tau_k}$. Note that at time $\tilde\tau_k$ the $\varepsilon_i$-coin lands head and, after arriving at $X_{\tau_k}$, the hyperplane of $X_{\tilde\tau_k}$ is never visited again. 
Conditioning on $X_{\tilde{\tau}_k}=x$, the law of $X_{\tau_k}$ is
$\mu_{\omega,\lambda,1}^x$, which is independent (under the environment measure $P$) of 
$\sigma(\omega(y,z):y\cdot e_1\le x\cdot e_1)$. \\
Moreover, after time $\tau_k$,
the path will never visit $\{y:y\cdot e_1\leq x\cdot e_1 +3L_1/4\}$. Therefore, 
$\tau_{k+1}-\tau_k$ is independent of what happened before $\tilde\tau_{k-1}$ and the inter-regeneration times form a $1$-dependent sequence.
\item
Since $(X_{\tilde\tau_{k+1}}-X_{\tau_k})_{k\ge 1}$ are iid and $(X_{\tau_{k+1}}-X_{\tilde\tau_{k+1}})\cdot e_1=1/\lambda_1$, the inter-regeneration distances
$\big((X_{\tau_{k+1}}-X_{\tau_k})\cdot e_1\big)_{k\ge 1}$ are iid.
\item 
From the construction, we see that a regeneration occurs after roughly a geometric number of levels. Thus we expect $(X_{\tau_{k+1}}-X_{\tau_k})\cdot e_1\sim c/\lambda$ and (by Lemma~\ref{hitting2}) $\tau_{k+1}-\tau_k\sim c/\lambda^2$.
\end{enumerate}

The above properties are made more precise in Lemma \ref{lemrenewal} below, the proof follows as in \cite{guo2012}. We introduce the $\sigma$-field
\[
\mathcal{G}_k
:=
\sigma\big(
\tilde{\tau}_k, (X_i)_{i\le\tilde{\tau}_k},(\omega(x,y))_{x\cdot e_1\le X_{\tilde{\tau}_k}\cdot e_1}
\big)
\]
and set 
\begin{equation}\label{e18}
	p_\lambda:=E \left[\sum_{y\in H_1}\mu_{\omega,\lambda,1}(y)P_{\omega,\lambda}^y(T_{-1/4}=\infty)\right].
\end{equation}
\begin{lem}\label{lemrenewal}
For any appropriate measurable sets $B_1, B_2$
and any event 
\[
B:=\{(X_i)_{i\ge 0}\in B_1, (\omega(x,y))_{x\cdot e_1>-L_1/4}\in B_2\},
\]
we have for $k\ge 1$,
\[
\bar{\mathbb P}_\lambda(B\circ\bar{\theta}_{\tau_k}|\mathcal{G}_k)
=
E\left[
\sum_{y\in H_1}
\mu_{\omega,\lambda,1}(y)
\bar{P}_{\omega,\lambda}^y(B\cap\{T_{-1/4}=\infty\})\right]
/p_\lambda.
\]
where  $\bar{\theta}_n$ is the time-shift defined by
\[
B\circ\bar{\theta}_n
=
\{(X_i)_{i\ge n}\in B_1, (\omega(x,y))_{(x-X_n)\cdot e_1>-L_0/\lambda_1}\in B_2\}.
\]
\end{lem}

We say that a sequence of random variables $(Y_i)_{i\in\mathbb{N}}$ is {\it $m$-dependent} ($m\in\mathbb{N}$) if 
\begin{equation*}
	\sigma(Y_i; 1\le i\le n)\mbox{ and }\sigma(Y_j; j> n+m)\mbox{ are independent,} \quad\forall n\in\mathbb{N}.
\end{equation*}
The law of large numbers and central limit theorem also hold for a stationary $m$-dependent sequence with finite means and variances, see \cite[Theorem 5.2]{Bi56}.
The following proposition is an immediate consequence of Lemma~\ref{lemrenewal}.

\begin{prop}\label{proprenewal}
Under $\bar{\mathbb P}_\lambda$,  $(X_{\tau_{n+1}}-X_{\tau_n})_{n\ge 1}$ and $(\tau_{n+1}-\tau_n)_{n\ge 1}$ are stationary 1-dependent sequences. Furthermore, for all $n\ge 1$, $(X_{\tau_{n+1}}-X_{\tau_n}, \tau_{n+1}-\tau_n)$ has the following law:
\begin{align*}
& \bar{\mathbb{P}}_\lambda (X_{\tau_{n+1}}-X_{\tau_n}\in\cdot,\tau_{n+1}-\tau_n\in\cdot)\\
&=
E\left[
\sum_{y\in H_1}
\mu_{\omega,\lambda,1}(y)
\bar{P}_{\omega,\lambda}^y(X_{\tau_1}\in\cdot,\tau_1\in\cdot,T_{-1/4}=\infty)\right]/p_\lambda.
\end{align*}
\end{prop}

\subsection{Moment estimates of regeneration times}
In this section we will show that the rescaled inter-regeneration times $\lambda^2(\tau_2-\tau_1)$ and inter-regeneration distances $\lambda(X_{\tau_2}-X_{\tau_1})$ have finite exponential moments.

\begin{thm}\label{distance_moment}
There exist constants $c_5,C_5>0$ such that for any $\omega$ and $\lambda\in(0,\lambda_0)$,
\[
\bar{E}_{\omega,\lambda}[\exp(c_5\lambda X_{\tau_1}\cdot e_1)]\leq C_5 < \infty\, .
\]
\end{thm}
\proof
First, observe that 
\begin{align*}
\bar{E}_{\omega,\lambda} [\exp(c\lambda X_{\tau_1}\cdot e_1)]
 &=\sum_{k\ge 1}\bar{E}_{\omega,\lambda}[\exp(c\lambda X_{S_k}\cdot e_1)\mathbbm{1}_{\{S_k<\infty,\tau_1=S_k\}}]\\
 &\le \sum_{k\ge 1}\bar{E}_{\omega,\lambda}[\exp(c\lambda X_{S_k}\cdot e_1)\mathbbm{1}_{\{S_k<\infty\}}].
\end{align*}
Next, by the definition of $S_k$, when $S_{k+1}<\infty$, we have (recall that $L_1=4L_0/\lambda_1$)
\[
X_{S_{1}}\cdot e_1=G_1 L_1
\]
and, recalling \eqref{e46},
\[
(X_{S_{k+1}}-X_{S_k})\cdot e_1=(N\circ\theta_{S_k}+G_k)  L_1
\quad\mbox{ for }k\ge 1,
\]
where $G_k:=\inf\{n\ge 1: \varepsilon_{n+M_k/L_1}=1\}$ is a geometric random variable with parameter $\beta$. Hence, taking $c>0$ small enough, we can achieve
\[
E[e^{c G_k}]\le \frac{9}{8}.
\]
 
Then, for $k\ge 1$, using the Markov property 
\begin{align*}
&\quad \bar{E}_{\omega,\lambda}[\exp(c X_{S_{k+1}}\cdot e_1/L_1)\mathbbm{1}_{\{S_{k+1}<\infty\}}]\\
&\le 
\frac{9}{8}\bar{E}_{\omega,\lambda}[\exp(cX_{S_k}\cdot e_1/L_1+c  N\circ\theta_{S_k})\mathbbm{1}_{\{S_k<\infty\}}]\\
&=\frac{9}{8}
\sum_x
\bar{E}_{\omega,\lambda}
\left[
\exp(c x\cdot e_1/L_1)\mathbbm{1}_{\{S_k<\infty,X_{S_k}=x\}}\right]
\bar{E}_{\omega,\lambda}^{x}
\left[ e^{cN}\mathbbm{1}_{\{T_{-1/4}<\infty\}}]
\right].
\end{align*}
Further,
for any $\omega$,
\begin{align*}
&\quad E_{\omega,\lambda}[e^{cN}\mathbbm{1}_{\{T_{-1/4}<\infty\}}]\\
&\le 
e^{c}P_{\omega,\lambda}(T_{-1/4}<\infty)+\sum_{n\ge 2}e^{cn}P_{\omega,\lambda}(N=n, T_{-1/4}< \infty)\\
&\leq
e^{c}P_{\omega,\lambda}(T_{-1/4}<\infty)+
\sum_{n\ge 2}\sum_z 
e^{cn}P_{\omega,\lambda}(X_{T_{n-1}}=z)P^z_{\omega,\lambda}(T_{-(n-1)-1/4}<T_{1})\\
& \stackrel{\mathclap{\rm{Cor.}\ \ref{cor1}}}{\le}\ 
e^{c}/2+\sum_{n\ge 1}(e^c/16)^n\le 7/8,
\end{align*}
where the last inequality is achieved by taking $c>0$ sufficiently small.
Thus, we conclude that 
taking $c>0$ sufficiently small, for $k\ge 1$,
\begin{align*}
\bar{E}_{\omega,\lambda}[\exp(cX_{S_{k+1}}\cdot e_1/L_1)\mathbbm{1}_{\{S_{k+1}<\infty\}}]
&\le 
\frac{63}{64}
\bar{E}_{\omega,\lambda}[\exp(c X_{S_k}\cdot e_1/L_1)\mathbbm 1_{S_k<\infty}].
\end{align*}
Therefore,
\begin{align*}
\bar{E}_{\omega,\lambda}[\exp(c X_{\tau_1}\cdot e_1/L_1)]
&\le 
\sum_{k\ge 1}\left(\frac{63}{64}\right)^{k-1}\bar{E}_{\omega,\lambda}[\exp(cX_{S_1}\cdot e_1/L_1)\mathbbm 1_{S_1<\infty}]\\
&=
64 E[e^{cG_1}]\le 72,
\end{align*}
which completes the proof. \qed

\begin{cor}\label{cor3}
There exist constants $c_6,C_6,c_7,C_7>0$ such that for any integer $n\ge 0$, any $\omega$ and $\lambda\in(0,\lambda_0)$,
\begin{equation}\label{e1}
\bar{\mb E}_\lambda[\exp(c_6\lambda(X_{\tau_{n+1}}-X_{\tau_n})\cdot e_1)]<C_6,
\end{equation}
\begin{equation}\label{e4}
\bar{\mb E}_\lambda\big[\exp\big(c_7\lambda^2(\tau_{n+1}-\tau_n)\big)\big]<C_7.
\end{equation}
\end{cor}
\proof 
To prove \eqref{e4}, it suffices to show that
\begin{equation}\label{e2}
\bar{\mb P}_\lambda(\tau_{n+1}-\tau_n>Cm/\lambda^2)\le Ce^{-c_6m}, \quad\forall m\geq 1, n\ge 0.
\end{equation}
If \eqref{e1} holds, we get the bound
\begin{align*}
\bar{\mb P}_\lambda(\tau_{n+1}-\tau_n>C_1m/\lambda^2)
&\le 
\bar{\mb P}_\lambda((X_{\tau_{n+1}}-X_{\tau_n})\cdot e_1\ge mL_1)
+\bar{\mathbb{E}}_\lambda \left[  P^{X_{\tau_n}}_{\omega,\lambda}(T_{m}>C_1m/\lambda^2)\right] \\
& \leq C_6e^{-cm} + e^{-c_1m}
\end{align*}
by Lemma~\ref{hitting2}. 
Thus we only need to prove \eqref{e1}. When $n=0$, inequality \eqref{e1} reduces to Theorem~\ref{distance_moment}. For $n\ge 1$, by Proposition~\ref{proprenewal},
\begin{align*}
\bar{\mb E}_\lambda[\exp(c_5\lambda(X_{\tau_{n+1}}-X_{\tau_n})\cdot e_1)]
\le 
2 E \left[
\sum_{y \in H_1}
\mu_{\omega,\lambda,1}(y) \bar E^y_{\omega,\lambda}[\exp( c_5\lambda X_{\tau_1}\cdot e_1) \mathbbm{1}_{\{T_{-1/4}=\infty\}}]
\right]
\le 2C_5,
\end{align*}
where we used again Theorem~\ref{distance_moment} and the fact (see Corollary~\ref{cor1}) that $P_{\omega,\lambda}(T_{-1/4}=\infty)\ge 1/2$ for all $\omega$ and $\lambda\in(0,\lambda_0)$. 
This proves \eqref{e1}. \qed\\

\begin{cor}\label{cor2}
Let $\lambda\in(0,\lambda_0)$.
\begin{enumerate}[(a)]
\item The speed $v(\lambda)$ satisfies $|v(\lambda)|\in(C\lambda, C'\lambda)$ for positive constants $C,C'$. 
\item The limit $Q_\lambda$ in \eqref{invmeasurelambda} exists $\mb P_\lambda$-almost surely for any bounded continuous $f$, and it defines an invariant measure for the process $(\bar\omega_n)_n$. Moreover,
{for any $f\in\mathcal F$, there exists $\lambda_0>0$ such that when $\lambda<\lambda_0$,}
\begin{equation}\label{e5}
Q_\lambda f=\bar{\mb E}_\lambda\left[\sum_{\tau_1\le i<\tau_2}
f(\bar\omega_i)
\right]\bigg/\bar{\mb E}_\lambda[\tau_2-\tau_1].
\end{equation}
\end{enumerate}
\end{cor}
\proof 
\begin{enumerate}[(a)]
\item Since the inter-regeneration distances and inter-regeneration times are stationary $1$-dependent sequences (Proposition~\ref{proprenewal}), and they have exponential moments (Corollary~\ref{cor3}), the law of large numbers gives
\[
v(\lambda) = \lim_{n\to\infty}\frac{X_n}{n}=\lim_{n\to\infty}\frac{X_{\tau_n}}{\tau_n}
=\frac{\bar{\mb E}_\lambda[X_{\tau_2}-X_{\tau_1}]}{\bar{\mb E}_\lambda[\tau_2-\tau_1]}
\]
Moreover, by Corollary~\ref{cor3} we have $|v(\lambda)|\ge \frac{1/\lambda}{c/\lambda^2}>C/\lambda$. On the other hand, Lemma~\ref{maxbound} implies
\begin{equation}\label{e8}
E_{\omega,\lambda}[T_{1/\lambda}] \geq \frac{s}{\lambda^2} P_{\omega,\lambda}(T_1\geq s/\lambda^2)
\ge \frac{s}{\lambda^2}P_{\omega,\lambda}(\max_{0\le t\le s/\lambda^2}|X_t|\le L_1)\ge s/2\lambda^2,
\end{equation}
for all $\omega\in\Omega$ and $\lambda\in(0,\lambda_0)$ and $s>0$ a sufficiently small constant.
Then we also have a lower bound for the moment of the inter-regeneration time
\begin{equation}\label{e7}
\bar{\mb E}_\lambda[\tau_{n+1}-\tau_n]\ge c/\lambda^2\quad\mbox{ for all }\lambda\in(0,\lambda_0), n\ge 0.
\end{equation}
So we have $|v(\lambda)|\le C'\lambda$.
\item
{
For a ballistic random walk in a uniformly elliptic finitely dependent random environment, recall that the regeneration times (which are different from the regeneration times in our paper) constructed by Comets and Zeitouni in \cite{cometszeitouni2004} divide both the path and the environment into i.i.d inter-regeneration pieces. This regeneration structure and the same argument as in the proof of Theorem~3.1 of \cite{sznizer1999} yields that the annealed law of $\bar\omega_n$ converges to an ergodic invariant measure, which we denote by $Q_\lambda$, of the sequence $(\bar\omega_n)$.
Hence
\[
Q_\lambda f=\lim_{n\to\infty}\bar{\mb E}_\lambda
\left[
\frac{1}{n}\sum_{i=0}^n f(\bar\omega_i) 
\right] \quad \text{for any $f\in\mathcal F$.}
\]
Suppose $f\in\mathcal F$ is  $\sigma(\omega(x,\cdot): |x|\le K)$-measurable for some $K>0$, then by Proposition~\ref{proprenewal}, when $\lambda<1/(4K)$, the sequence $\left(\sum_{\tau_k\le i<\tau_{k+1}}f(\bar\omega_i)\right), k\ge 1$, is 1-dependent and stationary. Therefore, by the moment estimates in Corollary~\ref{cor3} and the law of large numbers for 1-dependent stationary sequences, we have
\[
\lim_{n\to\infty}\bar{\mb E}_\lambda
\left[
\frac{1}{n}\sum_{i=0}^n f(\bar\omega_i)
\right]
=\lim_{k\to\infty}\bar{\mb E}_\lambda\left[\frac{1}{\tau_k}\sum_{i=0}^{\tau_k}f(\bar\omega_i)\right]
=\bar{\mb E}_\lambda\left[\sum_{\tau_1\le i<\tau_2}f(\bar\omega_i)\right]\bigg/\bar{\mb E}_\lambda[\tau_2-\tau_1].
\]
 \qed
}
\end{enumerate}


\section{Proof of Theorem \ref{newmain1}}\label{sec:gir}

Recall the notation $\mathbb E_{Q,\lambda}$ in Theorem~\ref{newmain1}.
We start by writing the expectation with respect to the unperturbed measure, 
\begin{align*}
\mathbb{E}_{Q,\lambda} \left[ \frac{\lambda}{t} \sum_{k=0}^{t/\lambda^2} f(\bar\omega_k)  \right] 
= \mathbb{E}_{Q,0} \left[ \frac{\lambda}{t} \sum_{k=0}^{t/\lambda^2} f(\bar\omega_k) \frac{dP_{\omega,\lambda}}{dP_\omega}(X_s; 0\leq s\leq t/\lambda^2)\right] 
\end{align*}
 and first study the Radon-Nikodym derivative 
\begin{equation}\label{radnikderiv}
G_\omega (\lambda, n) = \frac{\dd P_{\omega,\lambda}}{\dd P_\omega}(X_s; 0\leq s\leq n) .
\end{equation}


\subsection{An expansion of the Radon-Nikodym derivative}\label{girsanov1}

In this subsection we will derive a formula \eqref{densexp} for the density $G_\omega(\lambda, n)$. For a path $(x_0,\dots,x_n)$ with $x_0=0$ we have
\begin{align*}
\frac{P_{\omega,\lambda}((x_0,\dots,x_n))}{P_{\omega}((x_0,\dots,x_n))} &= \prod_{i=0}^{n-1} e^{\lambda \ell\cdot(x_{i+1}-x_i)} \frac{\sum_{|e|=1} \omega(x_i,x_i+e)}{\sum_{|e|=1} \omega(x_i,x_i+e)e^{\lambda\ell\cdot e}} \cr
&= 
\exp\left\{
	\lambda\ell\cdot x_n-\sum_{i=0}^{n-1}
	\log\left[\frac{\sum_{|e|=1}\omega(x_i,x_i+e)e^{\lambda\ell\cdot e}}{\sum_{|e|=1}\omega(x_i,x_i+e)}\right]
	\right\}.
\end{align*}

Note that for any $\omega\in\Omega$, 
\begin{align*}
\log\left[\frac{\sum_{|e|=1}\omega(0,e)e^{\lambda\ell\cdot e}}{\sum_{|e|=1}\omega(0,e)}\right]
&=\log\left[
\frac{\sum_{|e|=1}\omega(0,e)(1+\lambda\ell\cdot e+\frac{(\lambda\ell\cdot e)^2}{2}+o(\lambda^2))}{\sum_{|e|=1}\omega(0,e)}
\right]\\
&=\log\left[ 
1+\lambda E_\omega[X_1\cdot\ell]+\frac{\lambda^2}{2}E_\omega[(X_1\cdot\ell)^2]+o(\lambda^2)
\right]\\
&=\lambda E_\omega[X_1]\cdot\ell+\frac{\lambda^2}{2}{\rm Var}_\omega[X_1\cdot\ell]+o(\lambda^2),
\end{align*}
where in the last inequality we used the expansion $\log(1+x) = x-x^2/2+o(x^2)$.
Then, recalling $d(\omega,x)=E_\omega^x[X_1-X_0]$,
we obtain
\begin{align*}
\frac{P_{\omega,\lambda}((x_0,\dots,x_n))}{P_{\omega}((x_0,\dots,x_n))} 
&=
\exp\left\{ 
\lambda x_n\cdot\ell-\sum_{i=0}^{n-1} \left( \lambda d(\omega,x_i)\cdot\ell+\frac{\lambda^2}{2}{\rm Var}_{\theta_{x_i}\omega}[(X_1-X_0)\cdot\ell]+o(\lambda^2) \right) 
\right\}.
\end{align*}
Hence,  writing
\begin{equation} \label{defMMartingale}
M_n :=\bigg(X_{n} - \sum_{i=0}^{n-1} d(\omega,X_i)\bigg)\cdot\ell,
\qquad
D_\ell(\omega):={\rm Var}_\omega[(X_1-X_0)\cdot\ell],
\end{equation}
we conclude that
\begin{align}\label{densexp}
G(\lambda, n)
=\frac{\dd P_{\omega,\lambda}}{\dd P_\omega}(X_s; 0\le s\le n)
=\exp
\left\{ 
\lambda M_n
- \frac{\lambda^2}{2} \sum_{i=0}^{n-1}
D_\ell(\bar\omega_i)
 + n\cdot o(\lambda^2) 
\right\} .
\end{align}

\subsection{Weak convergence and Girsanov transform}\label{girsanov2}
In this subsection we will compute the limit of
\begin{align}\label{e17}
\mb E_{Q,\lambda}\left[ \lambda \sum_{k=0}^{t/\lambda^2}f(\bar\omega_k)\right] \nonumber
&=
\mb E_{Q,0}\left[ \lambda \sum_{k=0}^{t/\lambda^2}f(\bar\omega_k)G(\lambda,t/\lambda^2)\right] \nonumber\\
&=
\mb E_{Q,0}
\left[
\lambda \sum_{k=0}^{t/\lambda^2}f(\bar\omega_k)
\exp
\bigg\{
\lambda M_{t/\lambda^2}-\frac{\lambda^2}{2}\sum_{i=0}^{t/\lambda^2-1}D_\ell(\bar\omega_i)+t\cdot o(1)
\bigg\}
\right]
\end{align}
for $f\in\mathcal H_{-1}$ and any fixed $t>0$, as $\lambda\to 0$. First, we compute the limits of the terms in the expectation \eqref{e17}. 
Recall that 
for any $f\in \mathcal{H}_{-1}$, 
the process $\lambda\sum_{k=0}^{t/\lambda^2}f(\bar\omega_k)$ converges weakly (under $Q_0\times P_\omega$) to a Brownian motion $N_t^f$. 
Furthermore, notice that $M_n$ given in \eqref{defMMartingale} is a $P_\omega$-martingale whose increments are bounded and stationary with respect to $Q_0\times P_\omega$. 
Hence, the (joint) martingale CLT yields the joint convergence 
\begin{align} \label{jointweakconv}
\lambda \left(   
\sum_{k=0}^{t/\lambda^2} f(\bar\omega_k), M_{t/\lambda^2} 
\right)_{t \geq 0} 
\xrightarrow
[\lambda \rightarrow 0 ]{ } (N_t^f,N_t)_{t\geq 0} .
\end{align}
in distribution under $Q_0\times P_\omega$ to a 2-dimensional Brownian motion $(N_t^f, N_t)_{t\geq 0}$. 
Recall that the process $(\bar\omega_n)_{n\geq 0}$  is stationary and ergodic with respect to the initial measure $Q_0$ defined in \eqref{e16}. By the ergodic theorem, we have
\begin{align*}
\frac{\lambda^2}{2} \sum_{i=0}^{t/\lambda^2-1} 
D_\ell(\bar\omega_i)
\xrightarrow[\lambda \rightarrow 0 ]{} 
\frac{t}{2} E_{Q_0}\left[D_\ell(\omega) \right] 
\end{align*}
$Q_0\times P_\omega$-almost surely
and hence 
also $\mathbb{P}$-almost surely, where $E_{Q_0}$ denotes the expectation with respect to $Q_0$. 

Next, we will show that $G(\lambda,t/\lambda^2)$ is uniformly bounded in $L^p(P_{\omega})$, $p\ge 1$.
In fact, for any $p\ge 1$ and small enough $\lambda\le \lambda_0=\lambda_0(p,t)$,
\begin{equation}\label{e11}
E_{\omega}\left[ G(\lambda,t/\lambda^2)^p \right] \leq e^{p^2\frac{t}{2} + 1} .
\end{equation}
By the expansion of the Radon-Nikodym derivative in \eqref{densexp}, 
\[
\log G(p\lambda,t/\lambda^2)
-p\log G(\lambda,t/\lambda^2)
=
(p^2-p) \frac{\lambda^2}{2} \sum_{i=0}^{t/\lambda^2-1}  D_\ell(\bar\omega_i)+ C_{p,t} o(1) .
\]
Since $D_\ell(\omega)={\rm Var}_\omega[(X_1-X_0)\cdot\ell]\le 1$ for all $\omega$, we have
\[
|\log G(p\lambda,t/\lambda^2)
-p\log G(\lambda,t/\lambda^2)|
\le 
p^2 t/2+1
\]
for all $0<\lambda\le \lambda_0(p,t)$,  where $\lambda_0(p,t)$ is small enough. Inequality \eqref{e11} then follows by recalling that $E_\omega[G_\omega(p\lambda,t/\lambda^2)]=1$. 

Finally, we will show that
\begin{equation}\label{e19}
\lim_{\lambda\to 0} \mb E_{Q,\lambda}\left[ \lambda \sum_{k=0}^{t/\lambda^2}f(\bar\omega_k)\right]
=
t\mathrm{Cov}(N^f_1,N_1).
\end{equation}
Note that  the uniform integrability of $G(\lambda, t/\lambda^2)$ yields
\[
1=\mathbb E_{Q,0}[G(\lambda,t/\lambda^2)]
\xrightarrow[\lambda \rightarrow 0 ]{}
E\left[\exp\{{N_t}-\frac{t}{2}E_{Q_0}[D_\ell(\omega)]\}\right]
\] 
and hence we have necessarily $E[N_t^2]=t E_{Q_0}[D_\ell(\omega)]$. Furthermore, since $\lambda \sum_{k=0}^{t/\lambda^2} f(\bar\omega_k)$ is bounded in $L^2$ and the density $G (\lambda, t/\lambda^2)$ is bounded in $L^p$ for any $p \geq 1$, their product is uniformly integrable, which implies 
\begin{align*}
\lim_{\lambda\to 0}\mathbb{E}_{Q,\lambda} \left[ \lambda \sum_{k=0}^{t/\lambda^2} f(\bar\omega_k) \right]
&= \lim_{\lambda\to 0}
\mathbb{E}_{Q,0} \left[ \lambda\sum_{k=0}^{t/\lambda^2} f(\bar\omega_k) G(\lambda,t/\lambda^2) \right] \\
&\stackrel{\eqref{e17}}{=}
E\left[N_t^f e^{N_t -\frac{1}{2} E[N_t^2]} \right]
\end{align*}
By Girsanov's formula,
\begin{align*}
E\left[N_t^f e^{N_t -\frac{1}{2} E[N_t^2]} \right] = E\left[N_t^f + \mathrm{Cov}(N^f_t,N_t) \right] = \mathrm{Cov}(N^f_t,N_t),
\end{align*}
which proves \eqref{e19}.


\subsection{Remarks on the value of $\Lambda f$}

We have obtained in \eqref{e19} an expression for the operator $\Lambda$ in Theorem~\ref{newmain1}
\[
\Lambda f={\rm Cov}(N_1^f,N_1).
\]
We remark that this coincides with the definition of $\Lambda$ in \eqref{defLambda},
\[
\Lambda f=-{\rm Cov}(N_1^f,  N_1^d\cdot \ell),
\]
where $(N_t^d)_{t\ge 0}$ denotes the weak (Gaussian) limit of the process $\lambda\sum_{k=0}^{t/\lambda^2}d(\bar\omega_k,0)$ as $\lambda\to 0$. In the following, we denote by
\[
\Delta X_k:=X_{k+1}-X_k.
\] 
the increments of the random walk.
Noting that for any $n\ge 1$, the sequence \\
$(\Delta X_0,\ldots,\Delta X_{n-1}, \bar\omega_0,\ldots, \bar\omega_n)$ has the same distribution as 
$(-\Delta X_{n-1},\ldots,-\Delta X_0,\bar\omega_n,\ldots,\bar\omega_0)$
under  $Q_0\times P_\omega$, we have
\begin{equation} \label{orthogonality}
\mathbb E_{Q,0} \left[ (X_n-X_0)\cdot \left( \sum_{k=0}^n  f(\bar\omega_k) \right) \right] = \mathbb E_{Q,0} \left[ (X_0-X_n)\cdot \left( \sum_{k=0}^{n}  f(\bar\omega_k) \right) \right]=0.
\end{equation}
Consequently, 
\begin{align*}
\mathrm{Cov}(N^f_1,N_1) 
&= \lim_{n\to \infty } \frac{1}{n} 
{ \mathbb E_{Q,0}} \left[ \left( \sum_{k=0}^n f(\bar\omega_k)\right) \left( X_n - \sum_{i=0}^n  d(\omega,X_i) \right)\cdot\ell\right] \\
&= - \lim_{n\to \infty } \frac{1}{n} { \mathbb E_{Q,0}}  \left[  \left( \sum_{k=0}^n f(\bar\omega_k)\right) \left( \sum_{i=0}^n d(\omega,X_i) \cdot\ell\right)\right] \\
&= - {\rm Cov}(N^f_1, N^d_1\cdot\ell).
\end{align*}

We also remark that replacing the process $\lambda \sum_{k=0}^{t/\lambda^2}f(\bar\omega_k)$ by $\lambda X_{t/\lambda^2}$ in \eqref{e17},  the same argument gives
\begin{align}\label{e21}
\lim_{\lambda\to 0}\frac{1}{t}\mb E_{Q,\lambda}[\lambda X_{t/\lambda^2}]
&={\rm Cov}(B_1,N_1\cdot\ell)\nonumber\\
&=\lim_{n\to\infty}\frac{1}{n}\mb E_{Q,0}\left[ X_n\left( X_n - \sum_{i=0}^n  d(\omega,X_i) \right)\cdot\ell\right]\nonumber\\
&=\lim_{n\to\infty}\mb E_{Q,0}[X_n (X_n\cdot\ell)]=\Sigma\ell.
\end{align}

\section{A uniform LLN}\label{auLLN}

In this section, we show a quantitative result for the convergence to the steady state of a function $g$ of the environment seen from the particle and the increments of the process, provided that we can control maxima of the sum over $g$. In the next section, we will show that for functions $f\in \mathcal{F}$, we can control the maxima and can use the following theorem to prove Theorem \ref{newmain2}.

\begin{thm}\label{newmain2-2}
Let $g:\Omega\times\{e\in\mb Z^d:|e|\le 1\}\rightarrow \mathbb{R}$ be a function such that $g(\cdot, e)\in \mathcal F$ for any $|e|\le 1$. Assume that for all $\lambda$ smaller than some $\lambda_0>0$ and $n\ge 0$,
\begin{equation}\label{e20}
\bigg\lVert
\max_{0\leq m\le  n/\lambda^2}|\lambda \sum_{k=0}^m g(\bar\omega_k,\Delta X_k)|
\bigg\rVert_{L^{3/2}(Q_0\times P_{\omega,\lambda})}
\le  C e^{c\sqrt n}.
\end{equation}
Then, there exists a constant $C_8=C_8(\kappa,d,g)$ such that for any $t>1$ and $\lambda\in(0,\lambda_0)$,
\begin{equation}\label{newmain2f-2}
\Abs{
\frac{\frac{\lambda^2}{t}\mb E_{Q,\lambda}[\sum_{k=0}^{t/\lambda^2}g(\bar\omega_k,\Delta X_k)]- E_{\mb Q_\lambda} [g(\omega, \Delta X_0)]}{\lambda}
}
\le \frac{C_8}{t^{1/4}}.
\end{equation}
\end{thm}

Note that the sequence $(\bar\omega_k,\Delta X_k)_{k\ge 0}$ is stationary under the measure
\begin{equation}\label{notaQ}
\mathbb Q_\lambda:=Q_\lambda\times P_{\omega,\lambda}.
\end{equation}
Let us denote $\xi_n = \sum_{k=0}^n g(\bar{\omega}_k,\Delta X_k)$ and
\[  
\xi_{m,n}^*=\max_{m\le j\le n}\Abs{\sum_{k=m}^j g(\bar\omega_j,\Delta X_j)}, \quad \xi_n^*:=\xi_{0,n}^*.
\]
We first show the following consequence of inequality \eqref{e20}.


\begin{lem}\label{cor5}
Assume that \eqref{e20} holds.
 Then there exists a constant $C>0$ such that 
\[
\bar{\mathbb{E}}_{Q,\lambda}[(\lambda\xi^*_{\tau_2})^{4/3}]
\le C, \quad \forall \lambda\in(0,\lambda_0).
\]
\end{lem}

\proof
Let
$\bar{\mathbb{P}}_{Q,\lambda}
= Q_0 \times\bar{P}_{\omega,\lambda}$.
Since $(\xi^*_n)_{n\ge 0}$ is monotonically increasing, by Minkowski and H\"older's inequalities,
\begin{align*}
\norm{
\lambda \xi^*_{\tau_2}
}_{L^{4/3}(\bar{\mathbb{P}}_{Q,\lambda})}
&\le 
\left\lVert
\sum_{n=0}^\infty \lambda \xi^*_{n/\lambda^2,(n+1)/\lambda^2}
\mathbbm{1}_{\{\tau_2\ge n/\lambda^2\}}
\right\rVert_{L^{4/3}(\bar{\mathbb{P}}_{Q,\lambda})}\\
&\le 
\sum_{n=0}^\infty 
\left\lVert
\lambda\xi^*_{n/\lambda^2,(n+1)/\lambda^2}
\mathbbm{1}_{\{\tau_2\ge n/\lambda^2\}}
\right\rVert_{L^{4/3}(\bar{\mathbb{P}}_{Q,\lambda})}
\\
&\le
\sum_{n=0}^\infty
\left\lVert\lambda \xi^*_{n/\lambda^2,(n+1)/\lambda^2}
\right\rVert_{L^{3/2}(\bar{\mathbb{P}}_{Q,\lambda})}
\bar{\mathbb{P}}_{Q,\lambda}(\tau_2\ge n/\lambda^2)^{1/12}.
\end{align*}
By Corollary~\ref{cor3} and the fact that $\dd Q_0/\dd P$ is bounded (c.f. \eqref{e16}), we have
\[
\bar{\mathbb{P}}_{Q,\lambda}(\tau_2\ge n/\lambda^2)
\le 
C \bar{\mb P}_\lambda (\tau_2\ge n/\lambda^2)
\le 
C e^{-cn} \qquad\forall n\ge 0.
\]
Then, the lemma follows by applying inequality \eqref{e20}. 
\qed\\

An immediate consequence of \eqref{e7} and Lemma~\ref{cor5} is
\begin{equation}\label{e12}
{ |E_{\mb Q_\lambda}[g(\omega, \Delta X_0)]|}
\stackrel{\eqref{e5}}{=}
\Abs{
\frac{\bar{\mathbb{E}}_{\lambda}[\xi_{\tau_2}-\xi_{\tau_1}]}{\bar{\mathbb{E}}_{\lambda}[\tau_2-\tau_1]}
}
\le 
\frac{\bar{\mathbb{E}}_{Q,\lambda}[\xi^*_{\tau_1}+\xi^*_{\tau_2}]}{c/\lambda^2}
\le C \lambda.
\end{equation}

{
\noindent{\it Proof of Theorem~\ref{newmain2-2}}:\\
Note that by \eqref{e12}, inequality \eqref{e20} still holds if $g$ is replaced by $g-E_{\mb Q_\lambda}[g(\omega, \Delta X_0)]$. Thus without loss of generality, we assume that
\[
E_{\mb Q_\lambda}[g(\omega, \Delta X_0)]=0.
\]}
First, we will show that
\begin{equation}\label{e13}
\bar{\mathbb{E}}_{Q,\lambda}\left[
\max_{1\le k\le n}|\xi_{\tau_k}|
^{4/3}
\right]
\le 
C n, \qquad\forall n\ge 1.
\end{equation}
Indeed, for $m \ge 0$, let $N_m=\sum_{i=\tau_m}^{\tau_{m+1}-1}
g(\bar\omega_i,\Delta X_i)$. Then for any $m\ge 1$ by similar arguments as in the proof of \eqref{e5},
\[
\bar{\mb E}_\lambda[N_m]
= E_{\mb Q_\lambda}[g(\omega, \Delta X_0)] \bar{\mb E}_\lambda[\tau_2 - \tau_1] =0.
\]
Hence $(\sum_{m=0}^n N_m)_{n\ge 1}$ is a sequence with 1-dependent zero-mean increments under the measure $P\times P_{\omega,\lambda}$, and hence by \eqref{e16} also under $Q_0\times P_{\omega,\lambda}$. Moreover, by Lemma~\ref{cor5}, we have
\[
\bar{\mathbb{E}}_{Q,\lambda}[|\lambda N_m|^{4/3}]  
=
\bar{\mathbb{E}}_{Q,\lambda}[|\lambda N_1|^{4/3}] 
\leq 
2\bar{\mathbb{E}}_{Q,\lambda}[|\lambda \xi^*_{\tau_2}|^{4/3}] 
\leq C, \qquad\forall m\ge 0.
\]
By von Bahr-Esseen's inequality \cite[Theorem 2]{essen1965inequalities}, we have
\begin{align*}
\bar{\mathbb{E}}_{Q,\lambda}
\left[
\Abs{\sum_{\text{$m\le n$ is odd}}\lambda N_m}^{4/3}
\right]
\le 
2\sum_{\text{$m\le n$ is odd}}
\bar{\mathbb{E}}_{Q,\lambda}[|\lambda N_m|^{4/3}]
\le C n
\end{align*}
and then by Doob's $L^p$-martingale inequality, 
\[
\bar{\mathbb{E}}_{Q,\lambda}
\left[
\max_{k\le n}
\Abs{
\sum_{\text{$m\le k$ is odd}}\lambda N_m
}^{4/3}
\right]
\le C
\bar{\mathbb{E}}_{Q,\lambda}
\left[
\Abs{
\sum_{\text{$m\le n$ is odd}}\lambda N_m
}^{4/3}
\right]
\le 
C n.
\]
Similarly, we have
\[
\bar{\mathbb{E}}_{Q,\lambda}\left[
\max_{k\le n}
\Abs{
\sum_{\text{$m\le k$ is even}}\lambda N_m
}^{4/3}
\right]
\le 
C n.
\]
Combining these inequalities, we conclude that
\[
\bar{\mathbb{E}}_{Q,\lambda}\left[
\max_{k\le n}
\Abs{
\sum_{m=0}^k\lambda N_m
}^{4/3}
\right]
\le 
C n,
\]
which is exactly \eqref{e13}. 

Next, for $t>0$ fixed, we let $r=r(\lambda,t)\in\mathbb N$ be the integer that satisfies 
\[
\bar{\mb E}_\lambda[\tau_{r-1}]
\le 
4t/\lambda^2 
< 
\bar{\mb E}_\lambda[\tau_r].
\] 
We will show that
\begin{equation}\label{e23}
\bar{\mb P}_\lambda(t/\lambda^2\ge \tau_r)\le Ce^{-ct}.
\end{equation}

Note that by Corollary~\ref{cor3} and \eqref{e7},
we have 
\[
cr\le t<Cr.
\]
Since by Corollary~\ref{cor3}, the sequences $(\lambda^2\tau_1)_{\lambda\in(0,\lambda_0)}$ and $(\lambda^2(\tau_2-\tau_1))_{\lambda\in(0,\lambda_0)}$ are uniformly integrable (with respect to the measures $\bar{\mb P}_\lambda, \lambda\in(0,\lambda_0)$), by the lower bound \eqref{e7}, there exists a constant $M>0$ such that for all $n\ge 0$ and $\lambda\in(0,\lambda_0)$,
\[
\bar{\mb E}_\lambda[M\wedge\lambda^2(\tau_{n+1}-\tau_n)]\ge \frac{1}{2}\bar{\mb E}_\lambda[\lambda^2(\tau_{n+1}-\tau_n)].
\]
We set
\[
\bar\tau_n=\sum_{k=0}^{n-1}(\tau_{k+1}-\tau_k)\wedge (M/\lambda^2).
\]
Then, noting that $t/\lambda^2-\bar{\mb E}_\lambda[\bar\tau_r]\le t/\lambda^2-\frac{1}{2}\bar{\mb E}_\lambda[\tau_r]<-t/\lambda^2$, we have
\begin{align*}
\bar{\mb P}_\lambda(t/\lambda^2\ge \tau_r)
&\le 
\bar{\mb P}_\lambda(t/\lambda^2-\bar{\mb E}_\lambda[\bar\tau_r]\ge \bar\tau_r-\bar{\mb E}_\lambda[\bar\tau_r])\\
&\le 
\bar{\mb P}_\lambda(\lambda^2|\bar\tau_r-\bar{\mb E}_\lambda[\bar\tau_r]|\ge t)\\
&\le 
C\exp(-\frac{ct^2}{M^2 r})\le C e^{-ct},
\end{align*}
where in the third inequality we applied Azuma-Hoeffding's inequality to $\lambda^2(\bar\tau_r-\bar{\mb E}_\lambda[\bar\tau_r])$, which is a sum of bounded 1-dependent increments. Estimate \eqref{e23} is proved.

Now, by \eqref{e23},
\begin{align}\label{e14}
\bar{\mathbb{E}}_{Q,\lambda}\left[
 |\xi_{t/\lambda^2}|
\mathbbm{1}_{\{t/\lambda^2\ge \tau_r\}}\right]
\le 
\norm{\xi^*_{t/\lambda^2}}_{L^{3/2}(Q_0\times P_{\omega,\lambda})}
\bar{\mb P}_\lambda(t/\lambda^2\ge \tau_r)^{1/3}
\stackrel{\eqref{e20}}{\le} 
C/\lambda.
\end{align}

 On the other hand,
\begin{align*}
&\quad \bar{\mathbb{E}}_{Q,\lambda}
[
\xi_{t/\lambda^2}
\mathbbm{1}_{\{t/\lambda^2<\tau_r\}}
]\\
&\le 
\bar{\mathbb{E}}_{Q,\lambda}\left[
\max_{1\le k\le r}
|\xi_{\tau_k}|
\right]
+
\sum_{k=0}^{r-1}
\bar{\mathbb{E}}_{Q,\lambda}
\left[
\xi_{\tau_k,\tau_{k+1}}^*
\mathbbm{1}_{\{\tau_k\le t/\lambda^2<\tau_{k+1}\}}
\right]\\
&\stackrel{\mathclap{\eqref{e13}}}{\le}
\frac{C  r^{3/4}}{\lambda}
+
\sum_{k=0}^{r-1}
\norm{\xi_{\tau_k,\tau_{k+1}}^*}_{L^{4/3}(\bar{\mb P}_{Q,\lambda})}
\bar{\mathbb{P}}_{Q,\lambda}(\tau_k\le t/\lambda^2<\tau_{k+1})^{1/4}
\\
&\le
\frac{Ct^{3/4}}{\lambda}
+\frac{C}{\lambda}
\sum_{k=0}^{r-1}\bar{\mathbb{P}}_{Q,\lambda}(\tau_k\le t/\lambda^2<\tau_{k+1})^{1/4},
\end{align*}
where we used Lemma~\ref{cor5} and $r\leq Ct$ in the last step. Since by H\"{o}lder's inequality,
\[
\sum_{k=0}^{r-1}\bar{\mathbb{P}}_{Q,\lambda}(\tau_k\le t/\lambda^2<\tau_{k+1})^{1/4}
\le r^{3/4}\left(\sum_{k=0}^{r-1} \bar{\mathbb{P}}_{Q,\lambda}(\tau_k\le t/\lambda^2<\tau_{k+1})\right)^{1/4}
\le Ct^{3/4},
\]
we conclude that
\begin{equation}\label{e15}
\bar{\mathbb{E}}_{Q,\lambda}
\left[ |\xi_{t/\lambda^2}| \mathbbm{1}_{\{t/\lambda^2<\tau_r\}} \right]
\le 
C t^{3/4}/\lambda.
\end{equation}
Therefore, combining \eqref{e14} and \eqref{e15}, we arrive at
\begin{equation*}
 \mathbb E_{Q,\lambda}
\left[
|\xi_{t/\lambda^2}|
\right]
\le 
C t^{3/4}/\lambda, \qquad\forall t\ge 1.
\end{equation*}
This implies \eqref{newmain2f-2} and Theorem~\ref{newmain2-2} is proved.
\qed

\section{ Proofs of Theorem~\ref{newmain2} and the Einstein relation}\label{Proof of the Einstein relation}

\subsection{Proof of Theorem~\ref{newmain2}}
For $f\in\mathcal F$, let $g(\omega,x):=f(\omega)$. By Theorem~\ref{newmain2-2}, we only need to show that \eqref{e20} holds for $g$ in dimension $d\ge 3$. Note that now
$\xi_n = \sum_{k=0}^n f(\bar{\omega}_k)$ and recall that $\xi_{m,n}^*=\max_{m\le j\le n}|\sum_{k=m}^j f(\bar\omega_k)|$. In fact, we will obtain an estimate stronger than \eqref{e20}.
\begin{thm}\label{newmain2-3}
Let $d\ge 3$.
Let $f\in\mathcal F$ be a function that satisfies $Q_0 f=0$.
Then for all $\lambda\in [0,\lambda_0)$ and $n\ge 1$,
\[
\norm{\lambda\xi_{n/\lambda^2,(n+1)/\lambda^2}^*}
_{L^{3/2}(Q_0\times P_{\omega,\lambda})}
\le Cn^{6d+14}.
\]
\end{thm}

Our proof contains several steps.
\begin{enumerate}
\item First, we will obtain a moment estimate under the unperturbed measure
\begin{equation}\label{e9}
\norm{\xi^*_n}_{L^2(Q_0 \times P_{\omega})}\le C\sqrt n
\qquad\forall n\in\mb N.
\end{equation}

To prove inequality \eqref{e9}, recall the following theorem.
\begin{thm}[\cite{debuyermourrat2014, mourrat2011variance}]\label{thm1}
If $d\ge 3$ and $f:\Omega\to\mb R$ is an $L^2(Q_0)$ local function with $Q_0f=0$, then  for all $n\in\mathbb N$,
\[
E_{Q_0}[E_\omega[f(\bar\omega_n)]^2]\le C_f n^{-d/2}.
\]
Moreover, $(\sqrt n\xi_{nt})_{t\ge 0}$ converges weakly to a Brownian motion.
\end{thm}

Noting that $f(\bar\omega_n)$ is a stationary sequence under $Q_0 \times P_{\omega}$, by the maximal inequality for stationary sequences in \cite[Theorem 1]{peligrad2007maximal}, we have
\[
\norm{\xi^*_n}_{L^2(Q_0 \times P_{\omega})}\le C\sqrt n(\norm{f}_{L^2(Q_0)}+\delta_f),
\]
where $\delta_f=\sum_{m=1}^\infty m^{-3/2}\norm{E_\omega[\xi_m]}_{L^2(Q_0)}$. We only need to show that $\delta_f<\infty$. 
When $d\ge 3$, using the Cauchy-Schwarz inequality,
\begin{align*}
\norm{E_\omega[\xi_m]}_{L^2( Q_0)}^2
&=
E_{Q_0}\left[
\left(\sum_{j=1}^m E_\omega[f(\bar\omega_j)]\right)^2
\right]\\
&\le E_{Q_0}\left[
\sum_{j=1}^m \frac{1}{j^{d/4}} \cdot
\sum_{j=1}^m j^{d/4}(E_\omega[f(\bar\omega_j)])^2
\right]
\\
&\stackrel{\mathclap{{\rm Thm}~\ref{thm1}}}{\le}\ 
C_f 
\left(
\sum_{j=1}^m
j^{-d/4}
\right)^2\\
&\stackrel{\mathclap{d\ge 3}}{\le} C_f m^{1/2}.
\end{align*}
Hence $\delta_f\le C_f\sum_{m=1}^\infty m^{-3/2+1/4}<\infty$ for $d\ge 3$ and 
\eqref{e9} is proved. 

\item
Next, we will show that
\begin{equation}\label{e10}
\norm{\lambda\xi_{1/\lambda^2}^*}_{L^{5/3}(Q_0\times P_{\omega,\lambda})}\le C.
\end{equation}
Recalling the definition of the Radon-Nikodym derivative $G_\omega(\lambda,t)$ in 
\eqref{radnikderiv},
\begin{align*}
\norm{\lambda\xi^*_{1/\lambda^2}}_{L^{5/3}(Q_0\times P_{\omega,\lambda})}^{5/3}
&=
\mathbb E_{Q,0}
\left[
G(\lambda,1/\lambda^2)(\lambda\xi_{1/\lambda^2}^*)^{5/3}\right]
\\
&\le 
 \mathbb E_{Q,0} [G(\lambda,1/\lambda^2)^6]^{1/6}
 \mathbb E_{Q,0} [(\lambda\xi^{*}_{1/\lambda^2})^2]^{5/6}\\
&\stackrel{\mathclap{\eqref{e11},\, \eqref{e9}}}{\le} \quad 
C,
\end{align*}
which proves \eqref{e10}.

\item
For any fixed $\lambda\in(0,\lambda_0)$, let 
\begin{equation}\label{Def_h}
h(\omega)=h(\lambda,\omega)=E_{\omega,\lambda}[(\lambda\xi_{1/\lambda^2}^*)^{3/2}].
\end{equation}
Let $\square_n=\{x\in\mathbb Z^d: ||x||_{\infty} \le n\lfloor 1/\lambda\rfloor\}$ denote the box of side-length $2n\lfloor 1/\lambda\rfloor$ and set
\[
A_h(\omega):=\frac{1}{|\square_1|}\sum_{x\in\square_1}h(\theta_x\omega)
\]
to be an average of $h$ in the box $\square_1$.
We will show that for any integer $n\ge 0$,
\begin{equation}\label{box-average}
\mb E_{Q,\lambda}[h(\bar\omega_{n/\lambda^2})]
\le 
C\mb E_{Q,\lambda}[A_h(\bar\omega_{(n+1)/\lambda^2})
 +A_h(\bar\omega_{1+(n+1)/\lambda^2})
].
\end{equation}
To this end, we only need to prove that 
\[
\mb E_{Q,\lambda}[h(\bar\omega_{n/\lambda^2})]
\le 
C E_{Q,\lambda}[h(\theta_x \bar\omega_{(n+1)/\lambda^2})
+h(\theta_x \bar\omega_{1+(n+1)/\lambda^2})
], 
\quad \forall x\in \square_1.
\]
Further, noting that
\[
E_{\omega,\lambda}[h(\bar\omega_n)]
=
\sum_{y}P_{\omega,\lambda}(X_n=y)h(\theta_y\omega),
\]
it suffices to show that for any $y\in\mathbb Z^d$ and $x\in\square_1$,
\begin{equation}\label{parabolic-HE}
P_{\omega,\lambda}(X_{n/\lambda^2}=y)
\le C[P_{\omega,\lambda}(X_{(n+1)/\lambda^2}=y+x)
+P_{\omega,\lambda}(X_{1+(n+1)/\lambda^2}=y+x)].
\end{equation}
Indeed, for $k\in\mathbb Z$ and $y\in\mathbb Z^d$, set
\[
u(y,k)=P_{\omega,\lambda}^y(X_k=0).
\]
Then, by reversibility and uniform ellipticity,
\[
P_{\omega,\lambda}(X_k=y)
=
\frac{\sum_{e}\omega^\lambda(y,y+e)}{\sum_{e'}\omega^\lambda(0,e')}
u(y,k)
\asymp
e^{2\lambda y\cdot\ell} u(y,k),
\]
where $A\asymp B$ means $cB\le A\le CB$ for some constants $c,C>0$. Thus \eqref{parabolic-HE} is equivalent to 
\[
u(y,n/\lambda^2)\le C \left[u\left(y+x,(n+1)/\lambda^2\right)
+u\left(y+x,1+(n+1)/\lambda^2\right)
\right],
 \quad \forall x\in\square_1, y\in\mb Z^d,
\] 
which follows by (PHI) and
the fact that $u(\cdot,\cdot)$ satisfies the parabolic equation (PE) for the environment $\omega^\lambda$ in $B_{2\sqrt d/\lambda}(y)\times[(n-2)/\lambda^2,(n+2)/\lambda^2]$. 
Our proof of \eqref{box-average} is complete.

\item 
Note that by the Markov property,
\[
\norm{\lambda\xi_{n/\lambda^2,(n+1)/\lambda^2}^*}
_{L^{3/2}(Q_0\times P_{\omega,\lambda})}^{3/2}
=
\mb E_{Q,\lambda}[h(\bar\omega_{n/\lambda^2})]. 
\]
Thus by \eqref{box-average}, to prove Theorem~\ref{newmain2-3}, it suffices that $\mb E_{Q,\lambda} [A_h(\bar\omega_{n/\lambda^2})]\le C n^{9d+21}$, which by \eqref{e16} is equivalent to
\begin{equation}\label{e22}
\mb E_\lambda [A_h(\bar\omega_{n/\lambda^2})]\le C n^{9d+21}, \quad\forall n\in\mathbb N.
\end{equation}
\item 
We say that a box $\square_1(x):=x+\square_1$ centered at $x\in\mb Z^d$ is {\it $k$-good} (with respect to the environment $\omega$) if 
\[
A_h(\theta_z\omega)\le k^{9d+18} \qquad \text{for all }z\in\square_1(x).
\]
Otherwise, we say that $\square_1(x)$ is {\it $k$-bad}. We claim that 
\begin{equation}\label{bad}
P(\square_1 \text{ is $k$- bad})\le C k^{-10d-20}.
\end{equation}
Indeed, observing that $\sum_{y\in\square_2}h(\theta_y\omega)\big/|\square_2|\ge 2^{-d}A_h(\theta_z\omega)$ for all $z\in\square_1$, we have
\begin{align*}
P(\square_1 \text{ is $k$-bad})
&\le 
P(\sum_{y\in\square_2}h(\theta_y\omega)\big/|\square_2|\ge 2^{-d}k^{9d+18})\\
&\le C E\left[
\left(
\frac{1}{|\square_2|}\sum_{y\in\square_2}h(\theta_y\omega)
\right)^{10/9}
\right]\bigg/k^{10d+20}\\
&\le 
C E\left[
\frac{1}{|\square_2|}\sum_{y\in\square_2}h(\theta_y\omega)^{10/9}
\right]\bigg/k^{10d+20}\\
&=CE[h^{10/9}]/k^{10d+20}
\stackrel{\eqref{Def_h}}{\le} \mathbb E_\lambda[(\lambda\xi^*_{1/\lambda^2})^{5/3}]/k^{10d+20},
\end{align*}
where we used the translation-invariance of the measure $P$ in the second to last equality.
Display \eqref{bad} then follows by \eqref{e10} and the bounds of $\dd Q_0/\dd P$ in \eqref{e16}.  

\item 
Finally, we will prove \eqref{e22}. Clearly,
\begin{align*}
\mb E_\lambda[A_h(\bar\omega_{n/\lambda^2})]
&\le
\sum_{k=0}^\infty (k+1)^{9d+18}\mb P_\lambda\left(A_h(\bar\omega_{n/\lambda^2})\in[k^{9d+18}, (k+1)^{9d+18})\right)\\
&\le 1+C\sum_{k=1}^\infty
k^{9d+18}\mb P_\lambda(A_h(\bar\omega_{n/\lambda^2})\ge k^{9d+18}).
\end{align*}
Further, for each $k\in\mathbb N$, we can decompose the box $\square_k$ into $k^d$ boxes $(\square^{(i)})_{1\le i\le k^d}$ of side-length $2 \lfloor 1/\lambda\rfloor$. Hence
\begin{align*}
\mb P_\lambda(A_h(\bar\omega_{n/\lambda^2})\ge k^{9d+18})
&\le 
\mb P_\lambda(\max_{0\le s\le n/\lambda^2}|X_s|\ge k/\lambda)
+P(\text{one of the $k^d$ boxes $\square^{(i)}$ is $k$-bad})\\
&\stackrel{\text{Lem.}\ref{maxbound},\, \eqref{bad}}{\le }
C\frac{n^{9d+20}}{k^{9d+20}}+Ck^dk^{-10d-20}.
\end{align*}
Therefore, for any $n\in\mathbb N$,
\[
\mb E_\lambda[A_h(\bar\omega_{n/\lambda^2})]\le C\sum_{k=1}^\infty \frac{n^{9d+20}}{k^2}\le C n^{9d+20}.
\]
Inequality \eqref{e22} is proved.
\end{enumerate}
Our proof of Theorem~\ref{newmain2-3} is now complete.
\qed\\


\subsection{Proof of Theorem~\ref{ER}}
We note that by the ergodic theorem we can write the velocity as 
\begin{align*}
v(\lambda) = \lim_{n\to \infty} \frac{X_n}{n} =  Q_\lambda [ d(\omega^\lambda,0)]=E_{\mb Q_\lambda}[X_1], \qquad \mathbb{P}_\lambda\mbox{-a.s.}
\end{align*}
Let $g(\omega,e)=e\cdot e_i$ for a fixed unit vector $e_i$. Then, by Lemma~\ref{maxbound}, \eqref{e20} holds for this choice of $g$ and by Theorem~\ref{newmain2-2} we have
\[
\Abs{
\dfrac{\frac{\lambda^2}{t}\mb E_{Q,\lambda}[X_{t/\lambda^2}\cdot e_i]-E_{\mb Q_\lambda}[X_1\cdot e_i]}{\lambda}
}\le \frac{C}{t^{1/4}}.
\]
The collection of these inequalities for $i=1,\dots ,d$, together with \eqref{e21}, yields Theorem~\ref{ER}.

\subsection{Einstein relation as a corollary of Theorem~\ref{newER}}
We remark that Theorem~\ref{newER} can be considered a more general statement than the Einstein relation \eqref{ERconv}.
Indeed, since
\begin{align*}
v(\lambda)
=E_{Q_\lambda}E_{\omega,\lambda}[X_1]
&\stackrel{\eqref{densexp}}{=}
E_{Q_\lambda}E_\omega[e^{\lambda M_1+o(\lambda^2)}X_1]\\
&=
E_{Q_\lambda}E_\omega[(1+\lambda M_1)X_1]+o(\lambda)\\
&=
Q_\lambda[d(\omega,0)]+\lambda E_{Q_\lambda}E_\omega[M_1 X_1]+o(\lambda),
\end{align*}
we have by Theorem~\ref{newER}, applied to the collection of local functions $d(\omega,0)\cdot e_i$ for $i=1,\dots ,d$,
\[
\lim_{\lambda\to 0}\frac{v(\lambda)}{\lambda}
=
-{\rm Cov}(N_1^d, N_1^d\cdot\ell)+\mb E_{Q,0}[M_1 X_1].
\]
By the ergodic theorem and the fact that $\Delta X_k-d(\bar\omega_k,0)$ are $P_\omega$-martingale differences, 
\begin{align*}
{\rm Cov}(B_1-N_1^d, N_1)
&=
\lim_{n\to\infty}\frac{1}{n}\mb E_{Q,0}
\left[
\left (
X_n-\sum_{k=0}^{n-1}d(\bar\omega_k,0)
\right)
\left(
X_n-\sum_{k=0}^{n-1}d(\bar\omega_k,0)
\right)\cdot\ell
\right]\\
&=
\lim_{n\to\infty}\frac{1}{n}\mb E_{Q,0}
\left[
\sum_{k=0}^{n-1}(\Delta X_k-d(\bar\omega_k))(\Delta X_k-d(\bar\omega_k))\cdot\ell
\right]\\
&=
\mb E_{Q,0}[(X_1-d(\omega,0))(X_1-d(\omega,0))\cdot\ell]\\
&=\mb E_{Q,0}[X_1M_1].
\end{align*}
Therefore, 
\begin{align*}
\lim_{\lambda\to 0}\frac{v(\lambda)}{\lambda}
&=
-{\rm Cov}(N_1^d, N_1^d\cdot\ell)+{\rm Cov}(B_1-N_1^d, N_1)\\
&=-E[N_1^d(N_1^d\cdot\ell)]+E[(B_1-N_1^d)(B_1-N_1^d)\cdot\ell]
=E[B_1 (B_1\cdot\ell)]=\Sigma\ell ,
\end{align*}
since by \eqref{orthogonality}, $E[B_1N_1^d]=0$. Note  that we proved Theorem~\ref{newER} only for $d\ge 3$. Hence it does not cover the case $d=2$ of Theorem~\ref{ER}.


\section{Proof of the a-priori estimates}\label{Pape}

\subsection{Proof of Lemma \ref{hitting}}

We will prove Lemma~\ref{hitting} by contradiction. Let $u(x) = P_{\omega,\lambda}^x (T_{1,L} < T_{-1,L})$ and assume that for some $\omega$ we have 
\[
u(0)<\tfrac{2}{3}.
\] 
Recall that $T_{m,L}$ is the hitting time of $\{z\in \mathbb{Z}^d| z\cdot e_1 = nL/\lambda_1\}$. 
For a set $G\subset \mathbb{R}^d$, define its discrete {\it boundary} as
\[
\partial G=\{x\in\mb Z^d\cap G: x\sim y\text{ for some }y\in\mb Z^d\setminus G\}.
\] 
For a function $h:G\to \mathbb{R}$, denote its {\it Dirichlet energy} on $G$ as
\begin{align*}
\mathcal{E}(h,G) = \sum_{x,y\in G,\, x\sim y} \omega^\lambda(x,y) (h(x)-h(y))^2 .
\end{align*} 
We let $S_G=S_G(u,\omega)$ be the set of functions $v:G\to\mb R$ such that $v=u$ on $\partial G$. Since $u$ solves the elliptic equation (EE) in $\{x\in\mb Z^d: |x\cdot e_1|\le L/\lambda_1\}$, by Dirichlet's principle, we have
\[
\mathcal E(u,G)=\min_{v\in S_G}\mathcal E(v,G) \quad\text{ for any }G\subset\{x\in\mb Z^d: |x\cdot e_1|\le L/\lambda_1\}.
\]
We first find a lower bound for the Dirichlet energy of $u$ on the set $\Pi_0 = [-\tfrac{L}{\lambda_1},\tfrac{L}{\lambda_1}] \times [-\tfrac{1}{\lambda_1},\tfrac{1}{\lambda_1}]^{d-1} $ by ignoring all edges which are not in the direction of $e_1$. The conductance of such an edge in $\Pi_0$ connecting $(x_1,\bar{x})$ with $(x_1+1,\bar{x})$, where we write $x=(x_1,\bar{x})$ with $\bar{x}\in \mathbb{Z}^{d-1}$, is bounded from below by $c_\kappa e^{2\lambda \ell_1x_1}$. A lower bound for the energy is then
\begin{align*}
\mathcal{E}(u,\Pi_0) 
&\geq \sum_{||\bar{x}|| \leq 1/\lambda_1} \sum_{i=0}^{L/\lambda_1 -1} c_\kappa e^{2\lambda \ell_1i} (u(i+1,\bar{x})-u(i,\bar{x}))^2\\
&\ge C\sum_{||\bar{x}|| \leq 1/\lambda_1} \left( \sum_{i=0}^{L/\lambda_1 -1} u(i+1,\bar{x})-u(i,\bar{x})\right)^2 \left( \sum_{i=0}^{L/\lambda_1 -1} e^{-2\lambda \ell_1i} \right)^{-1} \\
&\ge C \sum_{||\bar{x}|| \leq 1/\lambda_1} (1-u(0,\bar{x}))^2 (1-e^{-2\lambda\ell_1}),
\end{align*}
where we used the Cauchy-Schwarz inequality in the second inequality. Note that $(1-e^{-2\lambda\ell_1})>c\lambda$ when $\lambda<\lambda_0$ for some small enough $\lambda_0>0$.
%
Since $u$ 
satisfies (EE) on $[-\tfrac{L}{\lambda_1},\tfrac{L}{\lambda_1}] \times [-\tfrac{2}{\lambda_1},\tfrac{2}{\lambda_1}]^{d-1} $, the elliptic Harnack inequality yields  $1-u(0,\bar{x})\geq C (1-u(0))$ whenever $|x|\leq 1/\lambda_1$. Thus
\begin{align}\label{lowerenergy}
\mathcal{E}(u,\Pi_0) \geq C \lambda \cdot \lambda^{-(d-1)} (1-u(0))^2 \geq  C\lambda^{2-d}
\end{align}
by our assumption.

Next, we will derive an upper bound for the energy. Define the (essentially one-dimensional) function
\begin{align*}
\bar{u}(x) = \frac{e^{2L} - e^{-2\lambda_1  x_1}}{e^{2L} - e^{-2L }} .
\end{align*}
Then $\bar{u}(x)=0$ when $x_1=-L/\lambda_1$ and $\bar{u}(x)=1$ when $x_1=L/\lambda_1$, that is, 
$\bar{u}$ satisfies in the strip $-L/\lambda_1 \leq e_1\cdot x \leq L/\lambda_1$ the same boundary conditions  as $u$. Let $R_1= \tfrac{L}{\lambda_1} \lfloor e^{L/d}\rfloor$ and $R_2=2R_1$. We will calculate energies on shifted sets $\Pi_1(y) = y+ [-\tfrac{L}{\lambda_1},\tfrac{L}{\lambda_1}] \times [-R_1,R_1]^{d-1} $ and  $\Pi_2(y) = y+ [-\tfrac{L}{\lambda_1},\tfrac{L}{\lambda_1}] \times [-R_2,R_2]^{d-1} $. For $L,\lambda$ fixed,  let
\begin{align*}
E_{L,\lambda} =\sup_{y\in\mathbb Z^d:y_1=0}\mathcal{E}(u,\Pi_1(y))<\infty
\end{align*}
and choose $y=y(L,\lambda)\in \mb Z^d$ such that with $\Pi_i=\Pi_i(y)$,
\begin{align*}
\mathcal{E}(u,\Pi_0)\leq  \mathcal{E}(u,\Pi_1(0))\le \mathcal{E}(u,\Pi_1) \quad \text{ and } \quad \mathcal{E}(u,\Pi_1) + \mathcal{E}(\bar{u},\Pi_1) > E_{L,\lambda} .
\end{align*}
We will show that for some positive constants $c_1,c_2$ independent of $L,\lambda$, 
\begin{align} \label{upperenergy}
\mathcal{E}(u,\Pi_1) \leq c_1 \lambda \cdot \lambda^{-(d-1)} e^{-c_2L} ,
\end{align} 
which contradicts \eqref{lowerenergy} if $L$ is large enough, since $\mathcal{E}(u,\Pi_0)\leq  \mathcal{E}(u,\Pi_1)$. 

To show \eqref{upperenergy}, set
\begin{align*}
v(x) = 
(1-d(x))\bar u(x)+d(x)u(x) \qquad\text{ for }x\in\Pi_2,
\end{align*}
where $d(x)=dist(x,\Pi_1)/ R_1 \le 1$. Note that $v=\bar u$ in $\Pi_1$ and $v=u$ on $\partial \Pi_2$.  By Dirichlet's principle,  $\mathcal{E}(u,\Pi_2)\leq \mathcal{E}(v,\Pi_2)$. 
For $x,y \in \Pi_2$ and $y\sim x$
\begin{align*}
v(x)-v(y) =(1-d(x))(\bar u(x)-\bar u(y))+ d(x)(u(x)-u(y))+(d(x)-d(y))(u(y)-\bar{u}(y)).
\end{align*}
{
Hence, observing $|d(x)-d(y)|\le R_1^{-1}\mathbbm{1}_{x,y\notin\Pi_1}$ for $x\sim y$, by Jensen's inequality,
\begin{align*}
(v(x)-v(y))^2 
\le& 
(1-d(x))(\bar{u}(x)-\bar{u}(y))^2+d(x)(u(x)-u(y))^2+R_1^{-2}\mathbbm{1}_{x,y\notin\Pi_1}(\bar{u}(y)-u(y))^2 \\
& + 2R_1^{-1}\mathbbm 1_{x,y\notin\Pi_1}|\bar{u}(y)-u(y)|\left[ (1-d(x))|\bar{u}(x)-\bar{u}(y)|+ d(x)|u(x)-u(y)|\right].
\end{align*}
Multiplying both sides by $\omega^\lambda(x,y)$ and summing over all $x,y\in \Pi_2$, this yields
\begin{align*}
\mathcal{E}(v,\Pi_2) & \leq \mathcal{E}(\bar{u},\Pi_2) + \mathcal{E}(u,\Pi_2\setminus \Pi_1^+) \\
& \tag{$a$} +  R_1^{-2} \sum_{y\in \Pi_2\setminus \Pi_1, x\sim y}  \omega^\lambda(x,y)(u(y)-\bar{u}(y))^2 \\
& \tag{$b$} + 2 R_1^{-1}  \sum_{y\in \Pi_2\setminus \Pi_1, x\sim y}  \omega^\lambda(x,y) |\bar{u}(y)-u(y)|\left[|\bar{u}(x)-\bar{u}(y)| + |u(x)-u(y)| \right],
\end{align*}
where $\Pi_1^+:=\Pi_1\setminus\partial\Pi_1$.
We will find upper bounds for the sums $(a)$ and $(b)$. }
Starting with the first one,
\begin{align*}
(a) \leq 2 R_1^{-2} \sum_{y\in \Pi_2\setminus \Pi_1, x\sim y}  \omega^\lambda(x,y)(u(y)-1)^2 + R_1^{-2} \sum_{y\in \Pi_2\setminus \Pi_1, x\sim y}  \omega^\lambda(x,y)(\bar{u}(y)-1)^2 
\end{align*}
with
\begin{align*}
(u(y)-1)^2 = \left( \sum_{i=y_1}^{L/\lambda_1} u(i,\bar{y}) - u(i+1,\bar{y}) \right)^2 \leq 2 \frac{L}{\lambda_1} \sum_{i=y_1}^{L/\lambda_1} (u(i,\bar{y}) - u(i+1,\bar{y}))^2
\end{align*}
by Cauchy-Schwarz's inequality.  Bounding $(\bar{u}(y)-1)^2$ analogously, we get
\begin{align*}
(a) & \leq 2  R_1^{-2} \frac{L}{\lambda_1} \left( \sum_{y\in \Pi_2\setminus \Pi_1, x\sim y}  \omega^\lambda(x,y) \sum_{i=y_1}^{L/\lambda_1} \big( u(i,\bar{y}) - u(i+1,\bar{y})\big)^2 +  \big(\bar u(i,\bar{y}) - \bar u(i+1,\bar{y})\big)^2 \right) \\
& \leq c  R_1^{-2} \frac{L}{\lambda_1} \left( \sum_{ y\in \Pi_2\setminus \Pi_1}   \sum_{i=y_1}^{L/\lambda_1} \omega^\lambda((i,\bar{y}),(i+1,\bar{y})) \left[ \big( u(i,\bar{y}) - u(i+1,\bar{y})\big)^2 +   \big(\bar u(i,\bar{y}) - \bar u(i+1,\bar{y})\big)^2\right] \right) \\
& \leq c  R_1^{-2} \left(\frac{L}{\lambda_1}\right)^2 \left( \mathcal{E}(u,\Pi_2\setminus \Pi_1^+) 
+ \mathcal{E}(\bar u,\Pi_2\setminus \Pi_1^+)\right) 
\end{align*}
For the summand $(b)$, H\"older's inequality yields
\begin{align*}
(b) & \leq 2 R_1^{-1} \left(  \sum_{y\in \Pi_2\setminus \Pi_1, x\sim y}  \omega^\lambda(x,y) (\bar{u}(y)-u(y))^2 \right)^{1/2} \\
& \qquad \qquad  \cdot \left(  \sum_{y\in \Pi_2\setminus \Pi_1, x\sim y}  \omega^\lambda(x,y) ({u}(x)-u(y))^2 + \sum_{y\in \Pi_2\setminus \Pi_1, x\sim y}  \omega^\lambda(x,y) (\bar{u}(x)-\bar{u}(y))^2 \right)^{1/2}
\end{align*}
The first sum can be bounded as we did for $(a)$ by $C(L/\lambda_1)^2(\mathcal{E}(u,\Pi_2\setminus \Pi_1^+) + \mathcal{E}(\bar u,\Pi_2\setminus \Pi_1^+))$ such that
\begin{align*}
(b) \leq c R_1^{-1} \frac{L}{\lambda_1} (\mathcal{E}(u,\Pi_2\setminus \Pi_1^+) + \mathcal{E}(\bar u,\Pi_2\setminus \Pi_1^+)) .
\end{align*}
Collecting the upper bounds for $(a),(b)$ and rearranging, using $ \mathcal{E}(u,\Pi_2)- \mathcal{E}(u,\Pi_2\setminus \Pi_1^+) =  \mathcal{E}(u,\Pi_1)${ and $\mathcal E(u, \Pi_2)\le \mathcal E(v,\Pi_2)$, we arrive at
\begin{equation}\label{e24}
\mathcal{E}(u,\Pi_1) 
\leq \mathcal{E}(\bar u,\Pi_2) 
+ c R_1^{-1}\frac{L}{\lambda_1} \left( \mathcal{E}(u,\Pi_2\setminus \Pi_1^+) + \mathcal{E}(\bar u,\Pi_2)\right).
\end{equation}
}
Next, we estimate the energy of $\bar{u}$. Since $\bar{u}(x)=\bar{u}(y)$ if $x_1=y_1$, we have
\begin{align*}
\mathcal{E}(\bar u,\Pi_2)&= \sum_{||\bar{x}||\leq R_2} \sum_{x_1=-L/\lambda_1}^{L/\lambda_1-1} \omega^\lambda((x_1,\bar{x}),(x_1+1,\bar{x})) 
\left( \frac{e^{2L} - e^{-2\lambda_1  x_1}}{e^{2L} - e^{-2L }} - \frac{e^{2L} - e^{-2\lambda_1  (x_1+1)}}{e^{2L} - e^{-2L }} \right)^2 \\
&\leq c R_2^{d-1} (e^{2L} - e^{-2L })^{-2} \sum_{x_1=-L/\lambda_1}^{L/\lambda_1-1} e^{2\lambda_1  x_1} \big( e^{-2\lambda_1  x_1} - e^{-2\lambda_1  (x_1+1)} \big)^2 ,
\end{align*}
where we used that $\omega^\lambda(x,y) \leq \kappa e^{\lambda \ell \cdot (x+y) } \leq \kappa e^{2\lambda_1  x_1}$. Now simple calculations give the upper bound (recall that $R_1= \tfrac{L}{\lambda_1} \lfloor e^{L/d}\rfloor$)
\begin{align*}
\mathcal{E}(\bar u,\Pi_2)& \leq c R_1^{d-1} (e^{2L} - e^{-2L })^{-2} \sum_{x_1=-L/\lambda_1}^{L/\lambda_1-1} e^{-2\lambda_1  x_1} \big( 1 - e^{-2\lambda_1 } \big)^2 \\
& \leq c R_1^{d-1} \lambda (e^{2L} - e^{-2L })^{-2} \sum_{x_1=-L/\lambda_1}^{L/\lambda_1-1} \lambda e^{-2\lambda_1  x_1}  \\
& \leq c R_1^{d-1} \lambda (e^{2L} - e^{-2L })^{-1} \\
& \leq c \left(\frac{L}{\lambda}\right)^{d-1}  \lambda e^{-L} .
\end{align*}
Finally, we use that $ \mathcal{E}(u,\Pi_2\setminus \Pi_1^+) \leq c E_{L,\lambda} \leq c (\mathcal{E}(\bar u,\Pi_1) + \mathcal{E}( u,\Pi_1))$ {and \eqref{e24} to obtain
\begin{align*}
\mathcal{E}( u,\Pi_1) 
 \leq 
\mathcal{E}(\bar u,\Pi_2) + ce^{-cL} (\mathcal{E}(\bar u,\Pi_2) + \mathcal{E}( u,\Pi_1)).
\end{align*} 
Therefore,
\begin{align*}
\mathcal{E}( u,\Pi_1)  
\le \dfrac{1+ce^{-cL}}{1-ce^{-cL}}\mathcal{E}(\bar u,\Pi_2) 
\leq 
C\left(\frac{L}{\lambda}\right)^{d-1}  \lambda e^{-L}.
\end{align*}
}
For $L$ large enough this implies \eqref{upperenergy} which then contradicts the lower bound \eqref{lowerenergy}.
\qed


\subsection{Proof of Lemma \ref{hitting2}}
\label{hitting2proof}

Let $\tilde{T}_1=T_1\wedge T_{-1}$. 
We will begin by estimating the lower bound of $P_{\omega,\lambda}(\tilde T_1>cn/\lambda^2)$ for all $n\in\mathbb N$ and $\lambda>0$ small enough. Note that this quantity depends only on the environments between the hyperplanes $\mathcal H_{-1}$ and $\mathcal H_1$, thus we let $a\in\Omega$ be an modified environment such that $a(x,y)=\omega^\lambda(x,y)$ for any $x,y\in\{z\in\mathbb Z^d: z\cdot e_1\in(-L_1,L_1)\}$ and $a$ satisfies $\frac{a(e)}{a(e')}<C$ for any bonds $e,e'$ in $\mb Z^d$.
Clearly,
\begin{align*}
P_{\omega,\lambda}\left( \tilde T_{1} \le \left(4L_1\right)^2 \right) 
=P_a\left( \tilde T_{1} \le \left(4L_1\right)^2 \right) 
\geq P_a (e_1\cdot X_{(4L_1)^2} >L_1).
\end{align*}
By the heat-kernel estimate in \cite[Theorem 3.1(i)]{grigor2002harnack},
there exist positive constants $\underline c_1,\underline c_2$ 
such that
\begin{align*}
P_a(X_{n}=y) \geq \frac{\underline c_1 m(y)}{V(0,\sqrt{n})} e^{-\underline c_2 |y|^2/n} ,
\end{align*}
whenever $|y|\leq n$ and $|y|_1+n$ is even, where  $m(y) = \sum_{z\sim y} a(y,z)$ and $V(x,r)= \sum_{z:|x-z|\leq r} m(z)$.
Define the set
\begin{align*}
A= \{ z \in \mathbb{Z}^d |\, \tfrac{4L_0}{\lambda_1} \leq e_1\cdot z \leq  \tfrac{8L_0}{\lambda_1} ,\ |e_i\cdot z | \leq \tfrac{4L_0}{\lambda_1} \text{ for } i=2,\ldots,d \} ,
\end{align*}
then $A\subset B_{4L_1}$ and $e_1\cdot z\geq L_1$ for all $z\in A$. Set $n=(4L_1)^2$ and let $A'$ be the set of $y\in A$ such that $|y|_1+(4L_1)^2$ is even, then we get
\begin{align*}
P_a\left( \tilde T_{1} \le \left(4L_1\right)^2 \right) 
& \geq \sum_{y\in A'} P_a (X_{(4L_1)^2}=y) \\
& \geq C\sum_{y\in A'} \frac{ m(y)}{V(0,4L_1)} e^{-\underline c_2 |y|^2/(16L/\lambda)^2}\\
& \geq C\sum_{y\in A'} \frac{e^{\lambda \ell \cdot y}}{\kappa^2 \sum_{z:|z|\leq 4L_1} e^{\lambda \ell \cdot z}} e^{-\underline c_2 } >C.
\end{align*}
This shows that there are positive constants $c,\delta>0$ such that
\begin{align} \label{proof:hitting1}
P_{\omega,\lambda}(\tilde T_{1} \le \tfrac{c}{\lambda^2} )>\delta \qquad\text{ for all }\omega\in\Omega.
\end{align}
Then, by the Markov property and  \eqref{proof:hitting1}
 we get for any $m\geq 1$,
\begin{align*} 
P_{\omega,\lambda}(\tilde T_{1} > \tfrac{mc}{\lambda^2} )\leq  \sup_{-L_1\leq x\cdot e_1 \leq L_1} P_{\omega,\lambda}^x(\tilde T_{1} >\tfrac{c}{\lambda^2})^m \leq (1-\delta)^m .
\end{align*}  
If we set $t_0=0$ and define recursively $t_{i+1}=t_i+\tilde{T}_1 \circ \theta_{t_i}$, where $\theta_n$ denotes the time shift of the trajectory (recall that $T_m$ is defined relative to the starting position), the exponential tail of $\tilde T_{1}$ implies by the exponential Markov inequality
\begin{align} \label{proof:hitting3}
P_{\omega,\lambda}(t_n > \tfrac{Cn}{\lambda^2} ) \leq e^{-2n} 
\end{align}
for some $C$ sufficiently large. We define the one-dimensional process $Y_n=(X_{t_n}\cdot e_1)/L_1$, which indicated the subsequent hyperplanes visited by $(X_n)_n$, then by Lemma \ref{hitting}, $Y_n$ jumps to the right with probability at least $\tfrac{2}{3}$ and then
\begin{align} \label{proof:hitting4}
P_{\omega,\lambda}\left(\sup_{k\leq Cn} Y_k < n\right) \leq e^{-2n} 
\end{align}
for $C$ sufficiently large. If we combine \eqref{proof:hitting3} and \eqref{proof:hitting4}, we obtain
\begin{align*}
P_{\omega,\lambda} ({T}_n > \tfrac{Cn}{\lambda^2}) \leq P_{\omega,\lambda}^x\left(\sup_{k\leq Cn} Y_k < n\right) + P_{\omega,\lambda}(t_n > \tfrac{Cn}{\lambda^2} ) \leq 2e^{-2n}.
\end{align*}
\qed


\subsection{Proof of Lemma \ref{maxbound}}
\label{maxboundproof}
 
We first find a uniform lower bound for the probability $P_{\omega,\lambda} (D_{1/\lambda} \geq \frac{r}{\lambda^2})$ for $r\leq 1$ small enough, where $D_{1/\lambda}$ is the exit time from the ball $B_{1/\lambda}$ of radius $\tfrac{1}{\lambda}$ in the 1-norm. 
Note that $P_{\omega,\lambda} (D_{1/\lambda} \geq \frac{r}{\lambda^2})$ depends only on the environments inside the box. Hence, as in subsection~\ref{hitting2proof}, we let $a\in\Omega$ be a modified environment such that $a(e)=\omega(e)$ for any bond $e$ in $B_{2/\lambda}$, and $a(e')/a(e'')<C$ for all bonds $e',e''$ in $\mb Z^d$. Then $P_{\omega,\lambda} ( D_{1/\lambda} \leq \tfrac{4r}{\lambda^2})=P_a( D_{1/\lambda} \leq \tfrac{4r}{\lambda^2})$, and
\begin{align*}
P_a ( D_{1/\lambda} \leq \tfrac{4r}{\lambda^2}) & 
\leq P_a (| X_{4r/\lambda^2}| \geq \tfrac{1}{2\lambda} ) + P_a(| X_{4r/\lambda^2}| < \tfrac{1}{2\lambda}, D_{1/\lambda} \leq \tfrac{4r}{\lambda^2}) \\
& \leq P_a (| X_{4r/\lambda^2}| \geq \tfrac{1}{2\lambda} ) + E_a \left[  P^{ X_{ D_{1/\lambda}}}_a (| X_{4r/\lambda^2- D_{1/\lambda}}| \geq \tfrac{1}{2\lambda})\, \mathbbm{1}_{\{ D_{1/\lambda} \leq \tfrac{4r}{\lambda^2}\} } \right] .
\end{align*}
By the heat kernel upper bound in \cite[Theorem 3.1(i)]{grigor2002harnack}, we get for $1/2\lambda_1 \leq k\leq 4r/\lambda_1^2$
\begin{align*}
P_a (|{X}_{k}|   \geq \tfrac{1}{2\lambda_1}) & = \sum_{m=1/2\lambda_1}^{k} \sum_{|x|=m} P_a ({X}_{k} =x)\\
& \leq C\sum_{m=1/2\lambda_1}^{k} \sum_{|x|=m} \frac{m(x)}{V(0,\sqrt{k})} e^{-\overline c_2 |x|^2/k} \\
& \leq C  k^{-d/2}   \sum_{m=1/2\lambda_1}^{k} m^{d-1}  e^{\lambda m-\overline c_2  m^2/k}\le C e^{-c/r},
\end{align*}
which is smaller than $\tfrac{1}{4}$ for $r$ small enough. This yields $P_{\omega,\lambda} ( D_{1/\lambda} \leq \tfrac{4r}{\lambda^2}) \leq \tfrac{1}{2}$ and so 
\begin{align*}
E_{\omega,\lambda}[e^{-\lambda^2D_{1/\lambda}}] \leq 1-\delta 
\end{align*}
 for some $\delta>0$ when $r,\lambda>0$ are small enough.
 Now we can proceed similarly to the proof of Lemma 4.5 in \cite{ganmatpia2012}:
\begin{align*}
E_{\omega,\lambda}\left[ \max_{0\leq s\leq t} |\lambda X_{s/\lambda^2}|^p \right] & = \int_0^\infty p y^{p-1} P_{\omega, \lambda} \left( \max_{0\leq s\leq t} |\lambda X_{s/\lambda^2}|^p \geq y \right) dy \\
& =  \int_0^\infty p y^{p-1} P_{\omega, \lambda} \left( D_{y/\lambda} \leq \tfrac{t}{\lambda^2} \right) dy \\
& \leq e^t \int_0^\infty p y^{p-1}E_{\omega,\lambda}\left[ e^{-\lambda^2D_{y/\lambda}}\right] dy 
\end{align*}
The exit time of the ball of radius $\tfrac{y}{\lambda^2}$ can be bounded as
\begin{align*}
D_{y/\lambda^2} \geq D_{1/\lambda^2} + D_{1/\lambda^2}\circ \theta_{D_{1/\lambda^2}} + \dots
+ D_{1/\lambda^2}\circ \theta_{D_{\lfloor y\rfloor /\lambda^2}} .
\end{align*}
The Markov property and the inequality $E_{\omega,\lambda}[e^{-\lambda^2D_{1/\lambda}}] \leq 1-\delta $ imply then for $t\leq 1$
\begin{align*}
E_{\omega,\lambda}\left[ \max_{0\leq s\leq t} |\lambda X_{s/\lambda^2}|^p \right] \leq 
e^t \int_0^\infty p y^{p-1}(1-\delta)^{\lfloor y \rfloor }  dy
\leq C , 
\end{align*}
with $C$ depending only on $p$, the bounds for the conductances and the dimension. 
This implies for $t\geq 1$
\begin{align*}
E_{\omega,\lambda}\left[ \max_{0\leq s\leq t} |\lambda X_{s/\lambda^2}|^p \right] \leq C \cdot t^p ,
\end{align*}
which is equivalent to the claimed inequality.
\qed

{\sl Acknowledgement:} The paper owes a lot to discussions with Pierre Mathieu. In particular, JN and XG thank Pierre for several invitations to Marseille, and for pointing us to the space $\mathcal{H}_{-1}$ and to reference \cite{mourrat2011variance}. Support of DFG (grant GA 582/8-1) is gratefully acknowledged.


\bigskip

\bibliographystyle{abbrv}
\bibliography{biber}

\begin{thebibliography}{}

\bibitem[Avena et~al., 2013]{lucarenatoflorian}
Avena, L., dos Santos, R.~S., and V{\"o}llering, F. (2013).
\newblock Transient random walk in symmetric exclusion: limit theorems and an
  {E}instein relation.
\newblock {\em ALEA Lat. Am. J. Probab. Math. Stat.}, 10(2):693--709.

\bibitem[Ben~Arous et~al., 2013]{ahoz2013einstein}
Ben~Arous, G., Hu, Y., Olla, S., and Zeitouni, O. (2013).
\newblock Einstein relation for biased random walk on {G}alton--{W}atson trees.
\newblock In {\em Annales de l'Institut Henri Poincar{\'e}, Probabilit{\'e}s et
  Statistiques}, volume~49, pages 698--721. Institut Henri Poincar{\'e}.

\bibitem[Billingsley, 1956]{Bi56}
Billingsley, P. (1956).
\newblock The invariance principle for dependent random variables.
\newblock {\em Trans. Amer. Math. Soc.}, 83:250--268.

\bibitem[Comets and Zeitouni, 2004]{cometszeitouni2004}
Comets, F. and Zeitouni, O. (2004).
\newblock A law of large numbers for random walks in random mixing
  environments.
\newblock {\em The Annals of Probability}, 32(1B):880--914.

\bibitem[de~Buyer and Mourrat, 2015]{debuyermourrat2014}
de~Buyer, P. and Mourrat, J.-C. (2015).
\newblock Diffusive decay of the environment viewed by the particle.
\newblock {\em Electron. Commun. Probab.}, 20:no. 23, 12.

\bibitem[De~Masi et~al., 1988]{masi1988}
De~Masi, A., Ferrari, P.~A., Goldstein, S., and Wick, W.~D. (1988).
\newblock An invariance principle for reversible {M}arkov processes.
  {A}pplications to random motions in random environments.
\newblock {\em Journal of Statistical Physics}, 55:787--855.

\bibitem[Delmotte, 1999]{delmotte1999}
Delmotte, T. (1999).
\newblock Parabolic {H}arnack inequality and estimates of {M}arkov chains on
  graphs.
\newblock {\em Rev. Mat. Iberoamericana}, 15(1):181--232.

\bibitem[Einstein, 1956]{einstein1956}
Einstein, A. (1956).
\newblock {\em Investigations on the theory of the {B}rownian movement.}
\newblock Edited with notes by R. F\"uhrt. Dower Publications, New York.

\bibitem[Gantert et~al., 2012]{ganmatpia2012}
Gantert, N., Mathieu, P., and Piatnitski, A. (2012).
\newblock Einstein relation for reversible diffusions in random environment.
\newblock {\em Comm. Pure. Appl. Math.}, 65:187--228.

\bibitem[Grigor'yan and Telcs, 2002]{grigor2002harnack}
Grigor'yan, A. and Telcs, A. (2002).
\newblock Harnack inequalities and sub-gaussian estimates for random walks.
\newblock {\em Mathematische Annalen}, 324(3):521--556.

\bibitem[Guo, 2012]{guo2012}
Guo, X. (2012).
\newblock Einstein relation for random walks in random environment.
\newblock {\em arXiv.org}, arXiv:1212.0255 [math.PR], to appear in: Ann. Prob.

\bibitem[Kipnis and Varadhan, 1986]{kipnis1986central}
Kipnis, C. and Varadhan, S. R.~S. (1986).
\newblock Central limit theorem for additive functionals of reversible {M}arkov
  processes and applications to simple exclusions.
\newblock {\em Communications in Mathematical Physics}, 104(1):1--19.

\bibitem[Komorowski et~al., 2012]{klo2012fluctuations}
Komorowski, T., Landim, C., and Olla, S. (2012).
\newblock {\em Fluctuations in Markov processes: time symmetry and martingale
  approximation}, volume 345.
\newblock Springer Science \& Business Media.

\bibitem[Komorowski and Olla, 2005a]{komorowskiolla2005einstein}
Komorowski, T. and Olla, S. (2005a).
\newblock Einstein relation for random walks in random environments.
\newblock {\em Stochastic {P}rocesses and their {A}pplications},
  115(8):1279--1301.

\bibitem[Komorowski and Olla, 2005b]{komorowskiolla2005mobility}
Komorowski, T. and Olla, S. (2005b).
\newblock On mobility and {E}instein relation for tracers in time-mixing random
  environments.
\newblock {\em Journal of {S}tatistical {P}hysics}, 118(3-4):407--435.

\bibitem[Lebowitz and Rost, 1994]{lebrost1994}
Lebowitz, J. and Rost, H. (1994).
\newblock The {E}instein relation for the displacement of a test particle in a
  random environment.
\newblock {\em Stochastic Process. Appl.}, 54:183--196.

\bibitem[Loulakis, 2002]{Loulakis}
Loulakis, M. (2002).
\newblock Einstein relation for a tagged particle in simple exclusion
  processes.
\newblock {\em Communications in mathematical physics}, 229(2):347--367.

\bibitem[Loulakis, 2005]{loulakis05}
Loulakis, M. (2005).
\newblock Mobility and {E}instein relation for a tagged particle in asymmetric
  mean zero random walk with simple exclusion.
\newblock {\em Ann. Inst. H. Poincar\'e Probab. Statist.}, 41(2):237--254.

\bibitem[Mourrat, 2011]{mourrat2011variance}
Mourrat, J.-C. (2011).
\newblock Variance decay for functionals of the environment viewed by the
  particle.
\newblock In {\em Annales de l'institut Henri Poincar{\'e} (B)}, volume~47,
  pages 294--327.

\bibitem[Peligrad et~al., 2007]{peligrad2007maximal}
Peligrad, M., Utev, S., and Wu, W. (2007).
\newblock A maximal $\mathbb{L}_p$-inequality for stationary sequences and its
  applications.
\newblock {\em Proceedings of the American Mathematical Society},
  135(2):541--550.

\bibitem[Shen, 2002]{shen2002}
Shen, L. (2002).
\newblock Asymptotic properties of certain anisotropic walks in random media.
\newblock {\em Annals of Applied Probability}, pages 477--510.

\bibitem[Spohn, 1989]{spohn}
Spohn, H. (1989).
\newblock Scaling limits for stochastic particle systems.
\newblock In {\em I{X}th {I}nternational {C}ongress on {M}athematical {P}hysics
  ({S}wansea, 1988)}, pages 272--275. Hilger, Bristol.

\bibitem[Sznitman and Zerner, 1999]{sznizer1999}
Sznitman, A.-S. and Zerner, M. (1999).
\newblock A law of large numbers for random walks in random environment.
\newblock {\em Ann. Prob.}, 27(4):1851--1869.

\bibitem[von Bahr and Esseen, 1965]{essen1965inequalities}
von Bahr, B. and Esseen, C.-G. (1965).
\newblock Inequalities for the {$r$}th absolute moment of a sum of random
  variables, {$1\leq r\leq 2$}.
\newblock {\em Ann. Math. Statist}, 36:299--303.

\end{thebibliography}

\bigskip
\bigskip

{\footnotesize

Nina Gantert: Technische Universit\"at M\"unchen,
Fakult\"at f\"ur Mathematik,
Boltzmannstra\ss e~3, 
85748~Garching bei M\"unchen,
Germany, gantert@ma.tum.de\\

Xiaoqin Guo: Purdue University,
Department of Mathematics,
150 N. University Street,
West Lafayette, IN 47907,
USA, guo297@purdue.edu \\

Jan Nagel: Technische Universit\"at M\"unchen,
Fakult\"at f\"ur Mathematik,
Boltzmannstra\ss e~3, 
85748~Garching bei M\"unchen,
Germany, jan.nagel@ma.tum.de\\
}

\end{document}